\documentclass[12pt]{article}

\usepackage{amsthm,amsfonts,amsmath,amssymb}
\usepackage[dvips]{graphicx}
\usepackage{alltt,makeidx,enumerate}
\allowdisplaybreaks
\newtheorem{thm}{Theorem}[section]
\newtheorem{lem}[thm]{Lemma}
\newtheorem{prop}[thm]{Proposition}

\theoremstyle{definition}

\newtheorem{remark}[thm]{Remark}

\numberwithin{equation}{section}
\newenvironment{prueba}{\noindent\textit{Proof}:}{  \qed\\\indent}

\newcommand{\dif}{\mathrm{d}}

\newcommand{\Div}{\mathrm{div\,}}

\newcommand{\Adj}{\mathrm{adj\,}}

\newcommand{\limite}[2]{\lim_{#1\rightarrow #2}}

\newcommand{\norm}[1]{\left\Vert#1\right\Vert}
\newcommand{\abs}[1]{\left\vert#1\right\vert}
\newcommand{\set}[1]{\left\{#1\right\}}
\newcommand{\mc}[1]{\mathcal{#1}}
\newcommand{\Real}{\mathbb R}
\newcommand{\eps}{\varepsilon}
\newcommand{\To}{\rightarrow}

\newcommand{\B}{\mathcal{B}}
\newcommand{\ts}[1]{\mathbf{#1}}
\newcommand{\tsg}[1]{\boldsymbol{#1}}

\newcommand{\pd}[2]{\frac{\partial#1}{\partial#2}}
\newcommand{\td}[2]{\frac{\dif#1}{\dif#2}}

\newcommand{\bx}{\mathbf{x}}
\newcommand{\bu}{\mathbf{u}}
\newcommand{\Rn}{\mathbb{R}^n}

\makeindex

\begin{document}

\title{A regularized penalty-multiplier method for approximating cavitation
solutions with prescribed cavity volume size}

\author{Pablo V. Negr\'on--Marrero\\
Department of Mathematics\\University of Puerto Rico\\
Humacao, PR 00791-4300\\ pablo.negron1@upr.edu\and
Jeyabal Sivaloganathan\\Department of Mathematical Sciences\\University of Bath,
Bath\\BA2
7AY, UK\\ masjs@bath.ac.uk}

\maketitle

\begin{abstract}
Let $\Omega\in\mathbb{R}^n$ be the region occupied by a body and let $\ts{x}_0$ be a flaw point in $\Omega$. Let $E(\cdot)$ be an energy functional (defined on some appropriate admissible set of
deformations of $\Omega$). For $V>0$ fixed, we let $\ts{u}_V$ be a minimizer of
$E(\cdot)$ among the set of deformations constrained to form a hole of volume
$V$ at $\ts{x}_0$. In this paper we describe a regularized penalty--multiplier method and its
convergence properties for the computation of both $\ts{u}_V$ and $E(\ts{u}_V)$.
In particular, we show 
that as the regularization parameter goes to zero, the regularized constrained
minimizers converge weakly in
$W^{1,p}(\Omega\setminus\overline{\B_{\delta}(\ts{x}_0)})$ to $\ts{u}_{V}$ for
any $\delta>0$. We describe as well the main features of a numerical scheme for
approximating $\ts{u}_V$ and $E(\ts{u}_V)$ and give a numerical example for the
case of a stored energy for an elastic fluid.
\end{abstract}


\section{Introduction}
When certain materials, such as rubber, are subjected to sufficiently
large triaxial loading, holes or bubbles begin to appear inside of the
stressed specimen (see, e.g., Gent and Lindley \cite{GeLi58}, Gent
\cite{Ge90}). The first variational model, based on the equations of nonlinear
elasticity, that predicts this phenomenon of void formation was given by Ball in
\cite{Ba82}. In this paper, Ball modelled a spherical body composed of an
isotropic ``soft" material and showed that, in the class of radial deformations,
any minimiser of the stored energy functional must open a hole at the
centre of the deformed ball for sufficiently large boundary displacements (the
phenomenon of cavitation). After this seminal paper, many others appeared on
different aspects of radial cavitation, e.g., \cite{AnNe87}, \cite{Ho92},
\cite{HoPo93}, \cite{HoPo95}, \cite{Si86a}, \cite{St85}, among
others. For results on general nonsymmetric cavitation in elasticity we refer to
\cite{JSp91}, \cite{MuSp95}, \cite{HMC} and \cite{SiSp2000a}.
A fundamental problem in studies of cavitation is to mathematically or
computationally characterize the boundary displacements for which
cavitation occurs. In \cite{NeSi2011b} the authors introduced the concept of the
\textit{volume derivative} as a tool for characterizing these boundary
displacements. For a large class of materials, the onset of cavitation--type
instabilities can be characterized as the zero level set of the volume
derivative. 

Central to any scheme for approximating the volume derivative is
the computation of minimizers of the corresponding stored energy functional
satisfying the constraint of forming or developing a hole inside
the body of a prescribed volume. This is the problem that concerns us in this
paper, in particular we describe a regularized penalty--multiplier
method and its convergence properties for the computation of such minimizers. This question was partially addressed in \cite{NeSi2011b} where a penalty method on a punctured domain is discussed. The convergence of this method as the penalty parameter increased to infinity was established, but no result was given as the radius of the spherical ``regularizing'' incision goes to zero. We provide such a result in this paper. Moreover, the use of a penalty--multiplier technique leads to a more a stable numerical scheme as compared to a standard penalty method. For more details on the penalty--multiplier method, also called \textit{augmented Lagrangians}, we refer to \cite{DL2008} and the references therein.

This paper is similar in spirit to \cite{SiSpTi2006} in which a regularization
method for cavitating solutions is described. However the problem considered in
\cite{SiSpTi2006} is without the constraint that the deformation opens a
hole inside the body of a prescribed volume. The presence of this nonlinear
volume constraint leads to various technical difficulties among which is that of
constructing variations satisfying this constraint. The use of
the penalty--multiplier technique eliminates the need for constructing such
variations by replacing the original constrained problem with a sequence of
unconstrained problems.

To introduce the results in the paper, consider a nonlinear hyperelastic
body occupying the bounded region $\Omega \subset \mathbb{R}^n$ in its
reference
state. A deformation of the body is a mapping $\bu:\Omega\rightarrow \Rn$
satisfying the local invertibility condition
\begin{equation}\label{in}
\det\nabla \bu (\bx ) > 0 \ \ \mathrm{a.e.} \ \bx \in \Omega.
\end{equation}
The energy stored in the deformed body under a deformation $\bu$ is
given by
\begin{equation}\label{energy}
E(\mathbf{u})=\int _{\Omega} W(\nabla \bu (\bx ) ) \ \mathrm{d}\bx,
\end{equation}
where $W: M^{n\times n}_+ \rightarrow \Real$ is the stored energy function of
the material and $M^{n\times n}_+$ denotes the set of $n\times n$ matrices
with positive determinant. In this paper we consider the
\textit{displacement boundary value problem} in which we fix a matrix
$\mathbf{A}\in M^{n\times n}_+$ and consider deformations satisfying
\begin{equation}\label{dispBC}
\mathbf{u}(\mathbf{x})=\mathbf{Ax} \mathrm{~for~
}\mathbf{x}\in\partial \Omega,
\end{equation}
We next fix a ``flaw" point $\mathbf{x}_0\in \Omega$, and for any fixed $V>0$  we take
the \textit{admissible set} of deformations to be
\begin{eqnarray}
\mathcal{A}_{\mathbf{A},V}&=&\big\{ \bu \in W^{1,p}(\Omega) \ | \
\mathrm{Det}\nabla
\bu = \mathrm{det}\nabla \bu \ \mathcal{L}^n+V \delta_{\mathbf{x}
_\mathbf{0}},\ \det\nabla \mathbf{u}>0 \ \mathrm{a.e.}, \nonumber \\
&&~~~~~\ \bu(\bx )=\mathbf{A}\bx \
\mathrm{on}\ \partial \Omega,\ \bu\ \mathrm{satisfies} \ \mathrm{INV}\mathrm{\
on \ }\Omega \big\}.
\label{admisset}
\end{eqnarray}
Here $\mathrm{Det}\nabla \bu$ denotes the distributional
determinant of $\bu$,
defined by
\begin{equation}\label{distjac}
<\mbox{Det}\nabla\ts{u},\phi>=-\dfrac{1}{n}\int_\Omega\nabla\phi
\cdot(\mbox{adj}\,\nabla\ts{u})\ts{u}\,\dif\ts{x},\quad\forall~~\phi\in
C^\infty_0(\Omega),
\end{equation}
$\mathcal{L}^n$ denotes $n$-dimensional Lebesgue measure, $p>n-1$,
$\delta_{\mathbf{x}_\mathbf{0}}$ denotes the Dirac measure
supported at $\mathbf{x}_\mathbf{0}\in \Omega$, and (INV)
denotes the  condition, relating to invertibility, introduced in
Definition 3.2 of \cite{MuSp95}. Results in \cite{SiSp2000a} give
conditions on the stored energy function $W$ under which a minimiser for
(\ref{energy}) exists on the set $\mathcal{A}_{\mathbf{A},V}$. The results of
Henao and Mora-Corral \cite{HMC} give conditions under which a minimiser also
exists in the case $p=n-1$ and their work in \cite{HMC2} includes justification
of the interpretation of $V$ in \eqref{admisset} as the volume of the hole
formed by the deformation. Hence, if $\bu\in\mathcal{A}_{\mathbf{A},V}$, then
the deformation $\bu$ produces a hole of volume $V$ in the deformed body. 

In Appendix \ref{regProb} we show that the requirement on deformations of
producing a hole of volume $V$ in the deformed body is equivalent to the
following \textit{integral constraint}: 
\[
\int_{\Omega}\det\nabla\ts{u}\,\dif\ts{x}=\left(\det\ts{A}\right)\,\abs{\Omega}
-V.
\]
(Here $\abs{\Omega}$ is the volume of $\Omega$.) Thus we replace the
minimization of \eqref{energy} over \eqref{admisset} with
\begin{equation}\label{AVmin0}
\left\{\begin{array}{ll}
\inf_{\ts{u}\in\mc{A}_{\ts{A}}}&E(\ts{u}),\\
\mbox{subject to 
}&c(\ts{u})=0,
\end{array}\right.
\end{equation}
where now
\begin{eqnarray*}
\mathcal{A}_{\mathbf{A}}&=&\big\{ \bu \in W^{1,p}(\Omega) \ | \
\exists\,\,\alpha\ge 0\mbox{ such that } \mathrm{Det}\nabla
\bu = \mathrm{det}\nabla \bu \ \mathcal{L}^n+\alpha \delta_{\mathbf{x}
_\mathbf{0}},\\
&&~~~~~\det\nabla \mathbf{u}>0 \ \mathrm{a.e.}, \ \bu(\bx )=\mathbf{A}\bx \
\mathrm{on}\ \partial \Omega,\ \bu\ \mathrm{satisfies} \ \mathrm{INV}\mathrm{\
on \ }\Omega \big\},
\end{eqnarray*}
and
\begin{equation}\label{const0}
c(\ts{u})=\int_{\Omega}\det\nabla\ts{u}\,\dif\ts{x}-
\left(\det\ts{A}\right)\,\abs{\Omega}+V.
\end{equation}
For any $\eps>0$, let
\[
\Omega_\eps=\Omega\setminus\overline{\B_{\eps}(\ts{x}_0)}.
\]
(Here and henceforth, we use the notation $\B_\eps(\ts{x}_0)$ for the open ball
of radius $\eps$ centered at $\ts{x}_0$.) The \textit{regularized constrained}
minimization problem is given by:
\begin{equation}\label{AVmineps}
\left\{\begin{array}{ll}
\inf_{\ts{u}\in\mc{A}^\eps_{\ts{A}}}&E_\eps(\ts{u}),\\
\mbox{subject to 
}&c_\eps(\ts{u})=0.
\end{array}\right.
\end{equation}
where
\begin{eqnarray}
E_\eps(\ts{u})&=&\int_{\Omega_\eps}W(\nabla\ts{u}(\ts{x}))\,\dif\ts{x},
\label{regenergy}\\
c_\eps(\ts{u})&=&\int_{\Omega_\eps}\det\nabla\ts{u}\,\dif\ts{x}-
\left(\det\ts{A}\right)\,\abs{\Omega}+V,\label{consteps}
\end{eqnarray}
and
\begin{eqnarray*}
\mathcal{A}^\eps_{\ts{A}} &=& \{ \ts{u} \in W^{1,p}(\Omega_\eps) \, |\,
\mathrm{Det}\nabla
\ts{u} = (\mathrm{det}\nabla \ts{u}) \, \mathcal{L}^n,\ \det \nabla \ts{u} >0 \
a.e., \\
&&~~~~~~~~~~\
\ts{u}(\ts{x} )=\ts{A}\ts{x} \ \mathrm{on}\ \partial\Omega,\ \ts{u}\
\mathrm{satisfies} \ \mathrm{INV} \},
\end{eqnarray*}
We set $\mc{A}^0_\ts{A}=\mc{A}_\ts{A}$ and $c_0=c$. In Appendix \ref{constAS}
we show that the admissible sets
\[
\mc{C}_{\ts{A}}^\eps\equiv\set{\ts{u}\in\mc{A}^\eps_{\ts{A}}\,|\,c_\eps(\ts{u}
)=0 },
\]
are nonempty for $\eps$ sufficiently small.
\begin{remark}
The hypotheses and results of \cite{SiSp2000a} are easily adapted to
prove that a (not necessarily unique) minimiser $\ts{u}_{V,\eps}$ of $E_\eps$ on
$\mathcal{A}^\eps_{\ts{A}}$ exists for each $\eps\ge0$ and $V>0$ small
enough.
\end{remark}

To compute approximations of the constrained problem \eqref{AVmineps}, we use a
penalty--multiplier method in which the functional in \eqref{regenergy} is
replaced with:
\begin{equation}
E_{\eps,\mu,\eta}(\ts{u})=E_\eps(\ts{u})+\mu\, c_\eps(\ts{u})
+\frac{1}{2}\eta\, c_\eps(\ts{u})^2.\label{Wpenalty}
\end{equation}
Here $\eta$ is a ``large'' positive parameter and $\mu\in\Real$.
Thus we replace the constrained problem \eqref{AVmineps} with the
``unconstrained'' problem:
\begin{equation}\label{AVminepseta}
\inf_{\ts{u}\in\mc{A}^\eps_{\ts{A}}}
~E_{\eps,\mu,\eta}(\ts { u } ).
\end{equation}
In Proposition \ref{prop:4.1} we show that a minimizer
$\ts{u}_{V,\eps,\mu,\eta}$ of \eqref{AVminepseta} actually exists. This result
is then used in Theorem \ref{thm:4.2} to show that for fixed $\eps,V>0$, there
exist sequences $\set{\mu_j}$, $\set{\eta_j}$ such that
$\set{\ts{u}_{V,\eps,\mu_j,\eta_j}}$ converges weakly in $W^{1,p}(\Omega_\eps)$
to a solution $\ts{u}_{V,\eps}$ of \eqref{AVmineps},
and $c_\eps(\ts{u}_{V,\eps,\mu_j,\eta_j})\To0$. We conclude Section
\ref{regpmProb} with a result on the weak form of the Euler--Lagrange
equations for the minimizer $\ts{u}_{V,\eps}$ (Theorem \ref{equivEQeps}) and a
result on the sensitivity of the minimum energy $E_\eps(\ts{u}_{V,\eps})$ to
variations on the boundary data $\ts{A}$ (Theorem \ref{sensieps}). 

In Section \ref{convgreg} we prove the main result of this paper, namely
Theorem \ref{rcmconv}. We show that for a sequence $\set{\eps_j}$ converging to
zero, the regularized constrained minimizers $\set{\ts{u}_{V,\eps_j}}$ converge
weakly in $W^{1,p}(\Omega_\delta)$ to a solution $\ts{u}_{V}$ of
\eqref{AVmin0}, for any $\delta>0$. The first part of the proof of this result,
dealing with the convergence and the existence of the limit $\ts{u}_{V}$, is
very similar to that in \cite[Theorem 4.1]{SiSpTi2006} and uses a
``diagonalization" argument. However showing that the limiting function
$\ts{u}_{V}$ is actually a solution of \eqref{AVmin0} is more subtle, again due
to the presence of the integral volume constraint in \eqref{AVmin0}. A key ingredient to get this result is Lemma \ref{lem:3.1} which shows that the energy of any function $\ts{u}\in\mc{C}_{\ts{A}}^0$ can be arbitrarily approximated by the energies of functions in $\mc{C}_{\ts{A}}^\eps$. The proof of this lemma involves the construction of certain diffeomorphisms of regularized sets where we needed to assume the convexity of $\Omega$. We close Section \ref{convgreg} with a result on the weak form of the Euler--Lagrange
equations for the minimizer $\ts{u}_{V}$ (Theorem \ref{equivEQ}).

Finally in Section \ref{numres} we describe the main features of a numerical
scheme for approximating solutions of the problem \eqref{AVmin0}. Also we give a
numerical example for the case of an elastic fluid ($\mu=0$ in
\eqref{sefexample}). For this class of materials and for a spherical domain, an
exact solution of \eqref{AVmin0} is known and thus we can check the various
convergence results in the paper in this case.

A simple class of polyconvex isotropic stored energy functions to which the
results in this paper can be applied is given by
\begin{equation}\label{sefexample}
W(\mathbf{F})=\frac \mu q\norm{\mathbf{F}}^q +h (\det \mathbf{F}),
\end{equation}
where $\kappa>0$, $q \in [n-1,n)$ and $h:(0,\infty )\rightarrow (0, \infty )$ is
such that
\begin{subequations}\label{hhyp}
\begin{eqnarray}
&h\mbox{ is a }C^2,\mbox{ convex function and}&\label{H1}\\
&h(\delta ) \rightarrow
\infty \mbox{ and } \frac{h(\delta )}{\delta } \rightarrow \infty\mbox{ as }
\delta \rightarrow 0, \infty\mbox{ respectively.}&\label{H2}
\end{eqnarray}
\end{subequations}
However, we note that the results of this paper apply to more general
polyconvex stored energy functions under varied hypotheses.

\section{The regularized penalty--multiplier method}\label{regpmProb}
In this section we study the unconstrained problem \eqref{AVminepseta}, in
particular the convergence properties of the penalty--multiplier method. We
assume that the stored energy function $W(\mathbf{F})$ satisfy the following:
\begin{enumerate}
\item[H1:]
(Polyconvexity) There exists with $G:(M_+^{n\times n})^{n-1}\times
(0,\infty )\rightarrow \mathbb{R}$ continuous and convex such that
\[
W(\mathbf{F})=\left\{\begin{array}{rcl} G(\mathbf{F}, \det\ts{F})&,&n=2,
\\G(\mathbf{F}, \Adj\mathbf{F}, \det\ts{F})&,&n=3.\end{array}\right.
\] 
\item[H2:]
(Growth) For $p\in (n-1,n)$, $c_1>0$, and a $C^2$ function $h$, we have that
\[
W(\mathbf{F})\geq K +c_1 |\mathbf{F}|^p +h(\det \mathbf{F})\quad
\mathrm{for}\quad \ts{F}\in M^{n\times n}_+,
\]
where the function $h$ satisfies conditions \eqref{hhyp}.
\end{enumerate}
We begin by showing that the minimizers in \eqref{AVminepseta} actually
exists.
\begin{prop}\label{prop:4.1}
For any $\mu\in\Real$, $\eta>0$, the infimum
\[
\inf_{\ts{u}\in\mc{A}^\eps_{\ts{A}}}
~E_{\eps,\mu,\eta}(\ts { u } ),
\]
exists and is attained for a function
$\ts{u}_{V,\eps,\mu,\eta}\in\mc{A}^\eps_{\ts{A}}$. Moreover, for any
$\delta>0$, the parameter $\eta$ can be chosen sufficiently large such that the
minimizer $\ts{u}_{V,\eps,\mu,\eta}$ satisfies that
$\abs{c_\eps(\ts{u}_{V,\eps,\mu,\eta})}<\delta$.
\end{prop}
\begin{prueba}
Since $\mc{A}^\eps_{\ts{A}}\ne\emptyset$, the infimum above is less than
$\infty$. If the infimum were $-\infty$, there would exists a sequence
$\set{\ts{u}_k}$ in $\mc{A}^\eps_{\ts{A}}$ such that $E_{\eps,\mu,\eta}(\ts { u
}_k )\To-\infty$ as $k\To\infty$, i.e.,
\[
\int_{\Omega_{\eps}}
W(\nabla\ts{u}_k(\ts{x}))\,\dif\ts{x}+\mu\, c_\eps(\ts{u}_k)
+\frac{1}{2}\eta\, c_\eps(\ts{u}_k)^2\To-\infty,\quad k\To\infty.
\]
Since the first and third terms above are bounded below, and for
$\ts{u}_k\in\mc{A}^\eps_{\ts{A}}$ we have that $c_\eps(\ts{u}_k)$ is bounded
below as well, it follows
that if $\mu\ge0$, the above limit can not be $-\infty$. If $\mu<0$, for the
above limit to be $-\infty$, we would
need to have (for a subsequence) that $c_\eps(\ts{u}_k)\To\infty$. In this last
case:
\begin{eqnarray*}
\int_{\Omega_{\eps}}
W(\nabla\ts{u}_k(\ts{x}))\,\dif\ts{x}+\mu\, c_\eps(\ts{u}_k)
+\frac{1}{2}\eta\, c_\eps(\ts{u}_k)^2&\ge&c_\eps(\ts{u}_k)\left(\mu 
+\frac{1}{2}\eta\, c_\eps(\ts{u}_k)\right),\\
&\To&\infty,\quad\mbox{as }k\To\infty,
\end{eqnarray*}
where we used that $\eta>0$. Thus in any case we arrive at a contradiction and
we must have that
\[
\inf_{\ts{u}\in\mc{A}^\eps_{\ts{A}}}
~E_{\eps,\mu,\eta}(\ts { u } )=g^*\in\Real.
\]

Let now $\set{\ts{u}_k}$ in $\mc{A}^\eps_{\ts{A}}$ be an infimizing sequence,
i.e., $E_{\eps,\mu,\eta}(\ts { u }_k )\To g^*$. An argument similar to the one
above implies that $\set{c_\eps(\ts{u}_k)}$ must be bounded. If $\mu
c_\eps(\ts{u}_k)\ge-L$ for all $k$, where $L>0$, then for $k$ sufficiently large
we
get that
\[
\int_{\Omega_{\eps}}
W(\nabla\ts{u}_k(\ts{x}))\,\dif\ts{x}-L\le g^*+1.
\]
It follows now from the growth hypotheses (H1)--(H2) that there exists a
subsequence $\set{\ts{u}_{k_j}}$ which converges weakly in
$W^{1,p}(\Omega_\eps)$ to a function $\ts{u}^*$, and that
$\set{\det\nabla\ts{u}_{k_j}}$ converges weakly in
$L^1(\Omega_\eps)$ to a function $\theta$. Since $p\in(n-1,n)$, it follows
from \cite[Theorem 4.2]{MuSp95}, that
$\ts{u}^*$ satisfies condition INV, $\theta=\det\nabla\ts{u}^*$,
and $\det\nabla\ts{u}^*>0$  almost everywhere. Thus
$\ts{u}^*\in\mc{A}^\eps_{\ts{A}}$.

Upon adapting the lower continuity results in \cite{SiSp2000a}, it follows that
$E_{\eps,\mu,\eta}$ is sequentially weakly lower semicontinuous. Thus we have
that
\[
E_{\eps,\mu,\eta}(\ts{u}^*)\le\liminf_{k_j}
E_{\eps,\mu,\eta}(\ts{u}_{k_j})=g^*=\inf_{\ts{u}\in\mc{A}^\eps_{\ts{A}}}
~E_{\eps,\mu,\eta}(\ts { u } ),
\]
i.e., that $\ts{u}_{V,\eps,\mu,\eta}\equiv\ts{u}^*\in\mc{A}^\eps_{\ts{A}}$ is a
minimizer.

For the last part of the proposition, we argue by contradiction. Thus we assume
that for some $\delta_0$ there exists a sequence $\eta_j\To\infty$ such that the
corresponding minimizers $\set{\ts{u}_j}$ satisfy
$\abs{c_\eps(\ts{u}_j)}\ge\delta_0$ for all $j$. Note that for all $j$,
\begin{equation}\label{energiesB}
E_{\eps,\mu,\eta_j}(\ts{u}_j)\le \left\{\begin{array}{ll}
\inf_{\ts{u}\in\mc{A}^\eps_{\ts{A}}}&\int_{\Omega_\eps}
W(\nabla\ts{u}(\ts{x}))\,\dif\ts{x},\\
\mbox{subject to 
}&c_\eps(\ts{u})=0.
\end{array}\right.\equiv f^*_\eps,
\end{equation}
Since $c_\eps(\ts{u}_j)\ge-C$ with $C>0$, if $\mu\ge0$, we get that
\[
f^*_\eps\ge
\mu c_\eps(\ts{u}_{j})+\frac{1}{2}\eta_{j}c_\eps(\ts{u}_{j})^2\ge -\mu C
+\frac{1}{2}\eta_{j}\delta_0^2\To\infty,
\]
as $j\To\infty$ which is a contradiction. If $\mu<0$ and
$\set{c_\eps(\ts{u}_j)}$ is
unbounded above, we can argue as before to get a contradiction as well. If
$\set{c_\eps(\ts{u}_j)}$ is bounded, then
\[
f^*_\eps\ge
\mu c_\eps(\ts{u}_{j})+\frac{1}{2}\eta_{j}c_\eps(\ts{u}_{j})^2\ge \mu
c_\eps(\ts{u}_{j})
+\frac{1}{2}\eta_{j}\delta_0^2\To\infty,
\]
which is a contradiction once again and this completes the proof.
\end{prueba}

We now show how to construct sequences $\set{\mu_j}$ and $\set{\eta_j}$ such
that the computed minimizers in \eqref{AVminepseta}, converge to the solution
of \eqref{AVmineps}.
\begin{thm}\label{thm:4.2}
Let the stored energy function $W$ satisfy the conditions H1--H2. Let
$\gamma\in(0,1)$, $\beta>1$, $\eta_1>0$, $\mu_1\in\Real$, and
$\ts{u}_0\in\mc{A}^\eps_{\ts{A}}$ be given. Let the sequences
$\set{\mu_j}$, $\set{\eta_j}$, and $\set{\ts{u}_j}$ be given by:
\begin{subequations}\label{penalITER}
\begin{eqnarray}
E_{\eps,\mu_j,\eta_j}(\ts{u}_{j})&=&\min_{\ts{u}\in\mc{A}^\eps_{\ts{A}}}
~E_{\eps,\mu_j,\eta_j}(\ts { u }
),\label{penalE}\\
\mu_{j+1}&=&\mu_{j}+\eta_{j}c_\eps(\ts{u}_{j}),\label{penalMU}\\
\eta_{j+1}&=& \left\{\begin{array}{rcl}
\eta_j&,&\mbox{if  }\abs{c_\eps(\ts{u}_j)}\le\gamma\abs{c_\eps(\ts{u}_{j-1})},\\
\beta\eta_j&,&\mbox{otherwise.}\end{array}\right.\label{penalETA}
\end{eqnarray}
\end{subequations}
Assume that $\set{\mu_j}$ is bounded. Then $c_\eps(\ts{u}_j)\To0$, and
$\set{\ts{u}_j}$ has a subsequence $\set{\ts{u}_{j_k}}$ that converges weakly in
$W^{1,p}(\Omega_\eps)$ to $\ts{u}_\eps\in\mc{A}^\eps_{\ts{A}}$ where
\begin{equation}\label{constP}
E_\eps(\ts{u}_{\eps})=\left\{\begin{array}{ll}
\min_{\ts{u}\in\mc{A}^\eps_{\ts{A}}}&\int_{\Omega_\eps}
W(\nabla\ts{u}(\ts{x}))\,\dif\ts{x},\\
\mathrm{subject~to 
}&c_\eps(\ts{u})=0.
\end{array}\right.
\end{equation}
\end{thm}
\begin{prueba}
By Proposition \ref{prop:4.1}, a function $\ts{u}_j\in\mc{A}^\eps_{\ts{A}}$
satisfying \eqref{penalE} exists for each $j$. From \eqref{energiesB} we get
that
\[
E_{\eps,\mu_j,\eta_j}(\ts{u}_{j})\le f^*_\eps,\quad\forall\,\,j.
\]
From this inequality and using that $W$ is nonnegative, we get that
\begin{equation}\label{auxPM2b}
\mu_jc_\eps(\ts{u}_j)+\frac{1}{2}\eta_jc_\eps(\ts{u}_j)^2\le
f^*_\eps,\quad\forall\,\,j.
\end{equation}
Note that the sequence $\set{\eta_j}$ is increasing. Thus in \eqref{penalETA} we
have two possibilities:
\begin{enumerate}[i)]
\item
the sequence $\set{\eta_j}$ remains bounded, in which case,
$\abs{c_\eps(\ts{u}_j)}\le\gamma\abs{c_\eps(\ts{u}_{j-1})}$ is satisfied for all
but but
finitely many indexes $j$. Clearly $c_\eps(\ts{u}_j)\To0$ in this case. 
\item
Otherwise (for a subsequence) $\eta_j\To\infty$, in which case \eqref{auxPM2b}
would imply that
$c_\eps(\ts{u}_j)\To0$.
\end{enumerate}
Thus in any case we have that $c_\eps(\ts{u}_j)\To0$.

If $\mu_j c_\eps(\ts{u}_j)\ge-L$ for all $j$, where $L>0$, then from
\eqref{energiesB} we get that
\[
\int_{\Omega_{\eps}}
W(\nabla\ts{u}_j(\ts{x}))\,\dif\ts{x}-L\le f^*_\eps.
\]
As in the proof of Proposition \ref{prop:4.1}, there exists a
subsequence
$\set{\ts{u}_{j_k}}$ which converges weakly in
$W^{1,p}(\Omega_\eps)$
to a function $\ts{u}_\eps$, that
$\set{\det\nabla\ts{u}_{j_k}}$ converges weakly in
$L^1(\Omega_\eps)$ to $\det\nabla\ts{u}_\eps$,  
$\ts{u}_\eps$ satisfies condition INV, and $\det\nabla\ts{u}_\eps>0$  almost
everywhere. Thus $\ts{u}_\eps\in\mc{A}^\eps_{\ts{A}}$ and
$c_\eps(\ts{u}_\eps)=0$.
Moreover, since $\mu_jc_\eps(\ts{u}_j)\To0$ by the assumed boundedness in
$\set{\mu_j}$, we have that
\[
f^*_\eps\le
E_\eps(\ts{u}_\eps)
\le\liminf_kE_{\eps,\mu_{j_k},
\eta_{j_k} }(\ts{u}_{j_k})\le f^*_\eps.
\]
It follows that 
\[
E_\eps(\ts{u}_{\eps})=\left\{\begin{array}{ll}
\min_{\ts{u}\in\mc{A}^\eps_{\ts{A}}}&\int_{\Omega_\eps}
W(\nabla\ts{u}(\ts{x}))\,\dif\ts{x},\\
\mathrm{subject~to 
}&c_\eps(\ts{u})=0.
\end{array}\right.
\]
\end{prueba}

Our next results give conditions under which the minimizer $\ts{u}_\eps$ in
\eqref{constP}, satisfies a weak form of the Euler-Lagrange equations for this
problem.
\begin{thm}\label{equivEQeps}
Let $\set{\ts{u}_j}$ be the sequence of Theorem \ref{thm:4.2} generated
according to \eqref{penalITER}, and $\set{\ts{u}_{j_k}}$ a subsequence that
converges weakly in $W^{1,p}(\Omega_\eps)$ to a solution $\ts{u}_\eps$ of
\eqref{constP}. Assume that there exist constants $K, \eps_0>0$ such that the
stored energy
function $W$ satisfies:
\begin{equation}\label{growthhyp2}
\abs{\td{W}{\ts{F}}(\mathbf{CF})\ts{F}^T}\leq K
 \left[
W(\mathbf{F})+1\right] \ \ \mathrm{for \ all \ }\mathbf{F}\in
\mathrm{M}^{n\times n}_+,
\end{equation}
whenever $\abs{\mathbf{C}-\mathbf{I}}<\eps_0$. Then $\set{\mu_j}$ has a
subsequence converging to $\mu_\eps$, where\footnote{$\mu_\eps$ is the
Lagrange multiplier corresponding to the volume constraint in \eqref{constP}
and is a measure of the Cauchy stress acting on the deformed inner cavity (cf.
\eqref{tracemult}).}
\begin{equation}\label{wfEL}
\int_{\Omega_\eps}\left[\td{W}{\ts{F}}(\nabla\ts{u}_\eps)
+\mu_\eps\left(\Adj\nabla\ts{u}_\eps\right)^T
\right]\cdot\nabla[\ts{v}(\ts{u}_\eps)]\,\dif\ts{x}=0,
\end{equation}
for all $\ts{v}\in C^1(\Real^n)$ with $\ts{v}=\ts{0}$ on
$\Real^n\setminus\mc{E}$, where $\mc{E}=\set{\ts{A}\ts{x}\,:\,\ts{x}\in
\Omega}$.
Moreover if $\ts{u}_\eps\in C^2(\Omega_\eps)\cap C^1(\overline{\Omega}_\eps)$
with
$\det\nabla\ts{u}_\eps>0$ in $\Omega_\eps$, then
\begin{subequations}\label{EQeqns}
\begin{eqnarray}
&\Div\left[\displaystyle\td{W}{\ts{F}}(\nabla\ts{u}_\eps)
+\mu_\eps\left(\Adj\nabla\ts{u}_\eps\right)^T\right]=\ts{0},\quad\mbox{in~~}
\Omega_\eps,&\label{EQreg}\\
&\ts{u}_\eps(\ts{x})=\ts{A}\ts{x}\mbox{~~on~~}\partial\Omega,&\label{EQbc1}\\
&\displaystyle\left[\td{W}{\ts{F}}(\nabla\ts{u}_\eps)
+\mu_\eps\left(\Adj\nabla\ts{u}_\eps\right)^T\right]\cdot\ts{n}=\ts{0}\mbox{
~~on~~}\partial\B_\eps(\ts{x}_0),&\label{EQbc2}\\
&\displaystyle\int_{\Omega_\eps}\det\nabla\ts{u}_\eps\,\dif\ts{x}=\frac{\omega_n
}{n}
\, \det\ts { A }
-V.&\label{EQconst}
\end{eqnarray}
\end{subequations}
\end{thm}
\begin{prueba}
To show \eqref{wfEL}, we first derive the corresponding equilibrium equation for
each $\ts{u}_j$. We use variations of $\ts{u}_j$ of the form 
$\ts{u}_s=\ts{u}_j+s\ts{v}(\ts{u}_j)$ where $\ts{v}\in C^1(\Real^n)$ with
$\ts{v}=\ts{0}$ on $\Real^n\setminus\mc{E}$. From \cite[Corollary
6.4]{SiSp2000a} it follows that for $s$ small enough, the function
$\ts{u}_s\in\mc{A}^\eps_{\ts{A}}$. (Note that the variation
$\ts{u}_s$ is not required to satisfy the constraint $c_\eps(\ts{u})=0$ as
$\ts{u}_j$ is a solution of an unconstrained problem.) To show \eqref{wfEL}
for $\ts{u}_j$, first note that
\begin{eqnarray*}
\int_{\Omega_\eps}[W(\nabla\ts{u}_s)\!\!\!&-&\!\!\!W(\nabla\ts{u}_j)]\dif\ts{x}
\\&=&
\int_{\Omega_\eps} \int_0^1\td{W}{t}(t\nabla\ts{u}_s+(1-t)\nabla\ts{u}_j)\,\dif
t\,\dif\ts{x},\\
&=&\int_{\Omega_\eps}
\int_0^1\td{W}{\ts{F}}(t\nabla\ts{u}_s+(1-t)\nabla\ts{u}_j)\cdot
(\nabla\ts{u}_s-\nabla\ts{u}_j)\,\dif t\,\dif\ts{x},\\
&=&s\int_{\Omega_\eps}\left[\int_0^1\td{W}{\ts{F}}([\ts{I}+st\nabla\ts{v}(\ts{
u}_j)]\nabla\ts{u}_j)\nabla\ts{u}_j^T\,\dif
t\right]\cdot\nabla\ts{v}(\ts{u}_j)\,\dif\ts{x}
\end{eqnarray*}
It follows now from \eqref{growthhyp2} that for $s$ small enough,
\[
\abs{\int_0^1\td{W}{\ts{F}}([\ts{I}+st\nabla\ts{v}(\ts{
u}_j))]\nabla\ts{u}_j)\nabla\ts{u}_j^T\,\dif
t}\le K[W(\nabla\ts{u}_j)+1]\in L^1(\Omega_\eps).
\]
Upon invoking the Dominated Convergence Theorem, we get that
\begin{equation}\label{aux121}
\limite{s}{0}\dfrac{1}{s}\int_{\Omega_\eps}\left[W(\nabla\ts{u}_s)-W(\nabla\ts{u
}
_j)\right]\dif\ts{x}=\int_{\Omega_\eps}\td{W}{\ts{F}}(\nabla\ts{u}_j)\nabla\ts{u
}
_j^T \cdot\nabla\ts{v}(\ts{u}_j)\, \dif\ts{x}
\end{equation}
Also
\begin{eqnarray*}
\mu_j[c_\eps(\ts{u}_s)-c_\eps(\ts{u}_j)]\!\!\!&+&\!\!\!\frac{1}{2}\eta_j[
c_\eps^2(\ts { u } _s)-
c_\eps^2(\ts{u}_j)]\\&=&[\mu_j+\frac{1}{2}\eta_j(c_\eps(\ts{u}_s)+c_\eps(\ts{u}
_j))][c_\eps(\ts{u}_s)-c_\eps(\ts{u}_j)].
\end{eqnarray*}
Now
\begin{eqnarray*}
c_\eps(\ts{u}_s)-c_\eps(\ts{u}_j)&=&\int_{\Omega_\eps}(\det\nabla\ts{u}
_s-\det\nabla\ts{u}_j)\,\dif\ts{x},\\
&=&\int_{\Omega_\eps}\int_0^1\td{}{t}\det(t\nabla\ts{u}_s+(1-t)\nabla\ts{u}_j)\,
\dif t\,\dif\ts{x},\\
&=&\int_{\Omega_\eps}\int_0^1[\Adj([\ts{I}+st\nabla\ts{v}(\ts{u}_j)]\nabla\ts{u}
_j)]^T\cdot(\nabla\ts{u}_s-\nabla\ts{u}_j)\,\dif t\,\dif\ts{x},\\
&=&s\int_{\Omega_\eps}\int_0^1[\Adj(\ts{I}+st\nabla\ts{v}(\ts{u}_j))]
^T(\Adj\nabla\ts { u }_j)^T\nabla\ts{u} _j^T\cdot\nabla\ts{v}(\ts{u}_j)\,\dif
t\,\dif\ts{x},\\
&=&s\int_{\Omega_\eps}\left[\int_0^1[\Adj(\ts{I}+st\nabla\ts{v}(\ts{u}
_j))]^T\,\dif t\right]\cdot\nabla\ts{v}(\ts{u}_j)\det\nabla\ts{u}_j\,\dif\ts{x}.
\end{eqnarray*}
It follows now since $\ts{v}\in C^1(\Real^n)$ with $\ts{v}=\ts{0}$ on
$\Real^n\setminus\mc{E}$, that
\begin{eqnarray*}
\limite{s}{0}\dfrac{1}{s}[c_\eps(\ts{u}_s)-c_\eps(\ts{u}_j)]
&=&\int_{\Omega_\eps}[\ts{I}
\cdot\nabla\ts{v}(\ts{u}
_j)]\det\nabla\ts{u}_j\,\dif\ts{x}.
\end{eqnarray*}
Combining this with \eqref{aux121} and using that $c_\eps(\ts{u}_s)\To
c_\eps(\ts{u}_j)$
as $s\To0$, we get that
\begin{eqnarray*}
\left.\td{}{s}E_{\eps,\mu_j,\eta_j}(\ts{u}_s)\right|_{s=0}&=&
\int_{\Omega_\eps}\bigg[\td{W}{\ts{F}}(\nabla\ts{u}_j)
\nabla\ts{u}_j^T\\&&+(\mu_j+\eta_jc_\eps(\ts{u}_j))(\det\nabla\ts{u}_j)\ts{I}
\bigg]\cdot\nabla\ts{v}(\ts{u}_j)\,\dif\ts{x}.
\end{eqnarray*}
Since $\ts{u}_j$ is a minimizer, we must have that
\begin{equation}\label{aux121a}
\int_{\Omega_\eps}\bigg[\td{W}{\ts{F}}(\nabla\ts{u}_j)
\nabla\ts{u}_j^T+(\mu_j+\eta_jc_\eps(\ts{u}_j))(\det\nabla\ts{u}_j)\ts{I}\bigg]
\cdot\nabla\ts { v } (\ts { u } _j)\,\dif\ts{x}=0,
\end{equation}
for all such $\ts{v}$'s. We now drop to the subsequence $\set{\ts{u}_{j_k}}$
that converges weakly in $W^{1,p}(\Omega_\eps)$ to $\ts{u}_\eps$ and with
$\det\nabla\ts{u}_{j_k}\rightharpoonup\det\nabla\ts{u}_\eps$ in
$L^1(\Omega_\eps)$.
Using \eqref{growthhyp2} and the Dominated Convergence Theorem once again, we
get that $\set{\mu_j}$ has a subsequence converging to $\mu_\eps$ where
\[
\int_{\Omega_\eps}\left[\td{W}{\ts{F}}(\nabla\ts{u}_\eps)
\nabla\ts{u}_\eps^T+\mu_\eps(\det\nabla\ts{u}_\eps)\ts{I}
\right]\cdot\nabla\ts{v}(\ts{u}_\eps)\,\dif\ts{x}=0.
\]
Since
$(\det\nabla\ts{u}_\eps)\ts{I}=(\Adj\nabla\ts{u}_\eps)^T\nabla\ts{u}_\eps^T$ and
$\nabla[\ts{v}(\ts{u}_\eps)]=\nabla\ts{v}(\ts{u}_\eps)\nabla\ts{u}_\eps$
, we get that the above equation can be written as
\[
\int_{\Omega_\eps}\left[\td{W}{\ts{F}}(\nabla\ts{u}_\eps)
+\mu_\eps(\Adj\nabla\ts{u}_\eps)^T
\right]\cdot\nabla[\ts{v}(\ts{u}_\eps)]\,\dif\ts{x}=0.
\]

Now assume that $\ts{u}_\eps\in C^2(\Omega_\eps)\cap
C^1(\overline{\Omega}_\eps)$
with
$\det\nabla\ts{u}_\eps>0$ in $\Omega_\eps$. Note that \eqref{EQbc1} and
\eqref{EQconst} follow from the fact that
$\ts{u}_\eps$ is a solution of \eqref{constP}. The proof that
\eqref{EQreg} holds is similar to the one given in \cite[Theorem 5.1]{SiSp2000a}
and thus we omit it. Now multiply \eqref{EQreg} by $\ts{v}(\ts{u}_\eps)$ where
$\ts{v}\in C^1(\Real^n)$ with $\ts{v}=\ts{0}$ on
$\Real^n\setminus\mc{E}$, and integrate by parts using \eqref{wfEL} to get that
\[
\int_{\partial\Omega_\eps}\left(\left[\td{W}{\ts{F}}(\nabla\ts{u}_\eps)
+\mu_\eps(\Adj\nabla\ts{u}_\eps)^T
\right]\cdot\ts{n}\right)\cdot\ts{v}(\ts{u}_\eps)\,\dif s(\ts{x})=0.
\]
Since the normal $\ts{n}$ to $\partial\Omega_\eps$ is mapped by $\ts{u}_\eps$ to
\[
\tilde{\ts{n}}(\ts{u}_\eps)=(\det\nabla\ts{u}_\eps)(\nabla\ts{u}_\eps)^{-T}
\ts{n},
\]
upon setting $\ts{y}=\ts{u}_\eps(\ts{x})$, the previous equation is equivalent
to:
\begin{equation}\label{aux119}
\int_{\ts{u}_\eps(\partial\Omega_\eps)}\left(\left[\ts{T}(\ts{y})+\mu_\eps\ts{I}
\right]\cdot\tilde{\ts{n}}(\ts{y})\right)\cdot\ts{v}(\ts{y})\,\dif s(\ts{y})=0,
\end{equation}
where the Cauchy stress tensor $\ts{T}(\ts{u}_\eps)$ is given by
\[
\ts{T}(\ts{u}_\eps)=(\det\nabla\ts{u}_\eps)^{-1}\td{W}{\ts{F}}(\nabla\ts{u}
_\eps)(\nabla\ts{u}_\eps)^T.
\]
Since $\ts{v}=\ts{0}$ on $\ts{u}_\eps(\partial\Omega)$, we get that
\eqref{aux119}
implies that
\[
\left[\ts{T}(\ts{y})+\mu_\eps\ts{I}
\right]\cdot\tilde{\ts{n}}(\ts{y})=\ts{0},\quad\forall\,\,
\ts{y}\in\ts{u}_\eps(\partial\B_\eps(\ts{x}_0)).
\]
which after changing variables back to $\Omega_\eps$ yields \eqref{EQbc2}.
\end{prueba}

We now study the sensitivity of the attained minimum value in \eqref{constP}
with respect to changes in the matrix $\ts{A}$. In the usual sensitivity
theorems of optimization theory, the parameters that change are in the right
hand sides of the constraints. In our problem however, the matrix
$\ts{A}$ appears both in the right hand side of the volume constraint and in the
displacement boundary condition on $\partial\Omega$ (c.f. \eqref{EQbc1},
\eqref{EQconst}). Thus our calculation picks up an additional term from
$\partial\Omega$. To emphasize the dependence of $\ts{u}_\eps$ on $\ts{A}$, we
use the notation $\ts{u}_\eps(\cdot,\ts{A})$.

\begin{thm}\label{sensieps}
Let $\ts{u}_\eps(\cdot,\ts{A})$ be the minimizer in \eqref{constP} and assume
that $\ts{u}_\eps(\cdot,\ts{A})\in C^2(\Omega_\eps)\cap
C^1(\overline{\Omega}_\eps)$
and
that $\ts{u}_\eps\in C^2(\Omega_\eps\times M_+^{n\times n})$. Then for
$\ts{A}=\mbox{diag}(\lambda_1,\ldots,\lambda_n)$, $\lambda_i>0$ for all $i$, we
have that
\begin{eqnarray}
\pd{}{\lambda_i}E_\eps(\ts{u}_\eps(\cdot,\ts{A}))&=&
\int_{\partial\Omega}x_i\ts{e}_i\cdot
\left[\td{W}{\ts{F}}(\nabla\ts{u}_\eps)
+\mu_\eps\left(\Adj\nabla\ts{u}_\eps\right)^T\right]\cdot\ts{n}\,\dif
s\nonumber\\&&-\mu_\eps\dfrac{\omega_n\det\ts{A}}{n\lambda_i},
\quad i=1,\ldots,n,\label{sensthm}
\end{eqnarray}
where $\set{\ts{e}_k}$ is the standard basis of $\Real^n$.
\end{thm}

\begin{prueba}
Let $\ts{u}_{\eps,i}=\pd{\ts{u}_\eps}{\lambda_i}$. By the assumed smoothness on
$\ts{u}_\eps$, we have that
\begin{equation}\label{aux213}
\pd{}{\lambda_i}E_\eps(\ts{u}_\eps(\cdot,\ts{A}))=\pd{}{\lambda_i}\int_{
\Omega_\eps}
 W(\nabla\ts{u}_\eps)\,\dif\ts{x}=\int_{\Omega_\eps}\td{W}{\ts{F}}(\nabla\ts{u}
_\eps)\cdot\nabla\ts{u}_{\eps,i}\,\dif\ts{x}.
\end{equation}
If we multiply \eqref{EQreg}
by $\ts{u}_{\eps,i}$, integrate by parts, and use the boundary condition
\eqref{EQbc2}, we get that
\begin{eqnarray*}
&\displaystyle\int_{\Omega_\eps}\left[\displaystyle\td{W}{\ts{F}}(\nabla\ts{u}
_\eps)
+\mu_\eps\left(\Adj\nabla\ts{u}_\eps\right)^T\right]\cdot\nabla\ts{u}_{\eps,i}\,
 \dif\ts{x}~~~~~~~~~~~~~&\\
&~~~~~~~~~~~~~~~~~~~~~~~~~~~=\displaystyle\int_{\partial\Omega}\left(\left[
\td { W } { \ts { F } } (\nabla\ts { u } _\eps)
+\mu_\eps\left(\Adj\nabla\ts{u}_\eps\right)^T\right]\cdot\ts{n}\right)\cdot
\ts{u}_{\eps,i}\,\dif s.&
\end{eqnarray*}
Since \eqref{EQbc1} implies that $\ts{u}_{\eps,i}(\ts{x})=x_i\ts{e}_i$ for
$\ts{x}\in\partial\Omega$, the above equation can be written as
\begin{eqnarray}
\displaystyle\int_{\Omega_\eps}\td{W}{\ts{F}}(\nabla\ts{u}_\eps)
\cdot\nabla\ts{u}_{\eps,i}\, \dif\ts{x}&=&
\displaystyle\int_{\partial\Omega}x_i\left(\left[
\td { W } { \ts { F } } (\nabla\ts { u } _\eps)
+\mu_\eps\left(\Adj\nabla\ts{u}_\eps\right)^T\right]\cdot\ts{n}\right)\cdot
\ts{e}_i\,\dif s\nonumber\\
&&-\mu_\eps\int_{\Omega_\eps}\left(\Adj\nabla\ts{u}_\eps\right)^T\cdot\nabla\ts
{u}_{\eps,i}\,\dif\ts{x}.\label{aux214}
\end{eqnarray}
We now differentiate \eqref{EQconst} with respect to $\lambda_i$ to get that
\[
\int_{\Omega_\eps}(\Adj\nabla\ts{u}_\eps)^T\cdot\nabla\ts{u}_{\eps,i}\,\dif\ts{x
}=
\dfrac{\omega_n}{n}\pd{}{\lambda_i}(\det\ts{A})=\dfrac{\omega_n\det\ts{A}}{
n\lambda_i}.
\]
Combining this with \eqref{aux213} and \eqref{aux214}, gives the result
\eqref{sensthm}.
\end{prueba}

\section{Convergence of the regularized constrained minimizers}\label{convgreg}
We now show that the regularized constrained minimizers whose existence is given
by Theorem \ref{thm:4.2}, converge to a solution of the ``non--regular"
constrained problem \eqref{AVmin0}. The first part of the proof of this result,
dealing with the convergence and the existence of the limit, is
very similar to that in \cite[Theorem 4.1]{SiSpTi2006} and consequently we
sketch most of it.  The second part in which we show that the limiting
function is actually a solution of \eqref{AVmin0} is more subtle,
again due to the treatment of the integral volume constraint in
\eqref{AVmin0}. In particular we need the following result which shows that
given any function $\ts{u}\in\mc{C}_{\ts{A}}^0$ and any sequence $\eps_j\To0$,
we can approximate approximate $\ts{u}$ with a sequence of functions
$\set{\hat{\ts{u}}_j}$ where $\hat{\ts{u}}_j\in\mc{C}_{\ts{A}}^{\eps_j}$ for all
$j$.
\begin{lem}\label{lem:3.1}
 Let $\Omega$ be a bounded, open, convex set, and let the stored energy function
$W$ satisfy the conditions H1--H2, and assume further that for some $K,\gamma_0>0$,
\begin{equation}\label{H3}
\abs{W(\ts{F}\ts{C})}\le K[W(\ts{F})+1],
\end{equation}
whenever $\norm{\ts{C}-\ts{I}}<\gamma_0$. Let $\ts{u}\in\mc{C}_{\ts{A}}^0$. Then
for any for any sequence $\set{\eps_j}$ with $\eps_j\To0$, there exists a
sequence of functions $\set{\hat{\ts{u}}_j}$ in $W^{1,p}(\Omega)$ such that
 \[
 \limite{j}{\infty}\int_{\Omega}W(\nabla\hat{\ts{u}}_j
(\ts{x}))\, \dif\ts{ x }
=\int_{\Omega}W(\nabla\ts{u}(\ts{x}))\,\dif\ts{x}.
\]
Moreover,
$\left.\hat{\ts{u}}_j\right|_{\Omega_{\eps_j}}\in\mc{C}_{\ts{A}}^{\eps_j}$.
\end{lem}
\begin{prueba}
 For any $0<\eta<1$ and $\eps>0$, we let
\[
\Omega_\eps^\eta=\set{\ts{x}\in\Omega_\eps\,:\,\mbox{dist}(x,\partial\Omega)>1-
\eta} .
\]
In terms of our previous notation, we have that $\Omega_\eps=\Omega_\eps^1$. Now
$\partial\Omega_\eps^\eta=\partial\B_\eps(\ts{x}_0)\cup\omega_\eta$ where
\[
\omega_\eta=\set{\ts{x}\in\Omega_\eps\,:\,\mbox{dist}(x,\partial\Omega)=1-
\eta}.
\]
For each $\ts{y}\in\omega_\eta$ there exists a unique
$\ts{x}(\ts{y})\in\partial \B_\eps(\ts{x}_0)$ such that
\[
\norm{\ts{x}(\ts{y})-\ts{y}}=\mbox{dist}(\ts{y},\partial\B_\eps(\ts{x}_0)).
\]
Since $\Omega$ is convex, the segment
\[
]\ts{x}(\ts{y}),\ts{y}[=\set{\gamma\ts{x}(\ts{y})+(1-\gamma)\ts{y}\,:\,
0<\gamma<1 },
\]
belongs to $\Omega_\eps^\eta$. Moreover, there exists a unique
$\ts{z}(\ts{y})\in \partial\Omega$ such that
\[
]\ts{x}(\ts{y}),\ts{y}[\subset]\ts{x}(\ts{y}),\ts{z}(\ts{y})[\subset\Omega_\eps.
\]
Since the segments $]\ts{x}(\ts{y}),\ts{y}[$, $]\ts{x}(\ts{y}),\ts{z}(\ts{y})[$
can be put into a one to one correspondence, after letting
$\ts{y}\in\omega_\eta$ to vary over $\omega_\eta$, this basically shows that we
can construct a diffeomorphism $\tsg{\Gamma}:\Omega_\eps\To\Omega_\eps^\eta$
such that $\det\nabla\tsg{\Gamma}>0$ over $\Omega_\eps$ and with $\tsg{\Gamma}(\ts{x})=\ts{x}$ for $\ts{x}\in\partial\B_\eps(\ts{x}_0)$. 
For any $\ts{u}\in\mc{C}_{\ts{A}}^0$, define $\hat{\ts{u}}_\eta$ on
$\Omega_\eps$
by:
\[
\hat{\ts{u}}_\eta(\ts{y})=\left\{\begin{array}{ll}
\ts{u}(\tsg{\Gamma}^{-1}(\ts{y})),&\ts{y}\in\Omega_\eps^\eta,\\
\ts{A}\ts{y},&\ts{y}\in\Omega_\eps\setminus\Omega_\eps^\eta.\end{array}\right.
\]
We now show that $\hat{\ts{u}}_\eta\in\mc{C}_{\ts{A}}^{\eps}$ for $\eta$
(depending on $\eps$) sufficiently close to 1.
\begin{enumerate}[i)]
\item
Clearly $\hat{\ts{u}}_\eta(\ts{x})=\ts{A}\ts{x}$ for $\ts{x}\in\partial\Omega$,
and
\[
\det\left[\nabla(\ts{u}(\tsg{\Gamma}^{-1}(\ts{y})))\right]= \det\nabla
\ts{u}(\tsg{\Gamma}^{-1}(\ts{y}))\,\det\nabla\tsg{\Gamma}^{-1}(\ts{y})>0,
\]
a.e. on $\Omega_\eps$ since $\det\nabla\ts{u}>0$ a.e. on $\Omega$, and
$\det\nabla\tsg{\Gamma}>0$ implies that $\det\nabla\tsg{\Gamma}^{-1}>0$ on
$\Omega_\eps$. Clearly $\hat{\ts{u}}_\eta$ satisfies INV in $\Omega_\eps$ and since $\delta_{\ts{x}_0}(\Omega_\eps)=0$, it follows that $\hat{\ts{u}}_\eta\in\mc{A}^\eps_{\ts{A}}$.
\item
Since $\ts{u}\in\mc{C}_{\ts{A}}^0$, we have that
\[
\int_{\Omega}\det\nabla\ts{u}\,\dif\ts{x}=\left(\det\ts{A}\right)\,\abs{\Omega}
-V.
\]
Hence
\[
\int_{\Omega_\eps}\det\nabla\ts{u}\,\dif\ts{x}=\left(\det\ts{A}\right)\,\abs{
\Omega}
-V_\eps,\quad
V_\eps=V+\int_{\B_\eps(\ts{x}_0)}\det\nabla\ts{u}\,\dif\ts{x}.
\]
Note that $V_\eps>V$, and $V_\eps\searrow V$ as $\eps\searrow0$. Now
\begin{eqnarray*}
\int_{\Omega_\eps}\det\nabla\hat{\ts{u}}_\eta\,\dif\ts{\ts{y}}&=&
\int_{\Omega_\eps^\eta}\det\nabla\ts{u}(\tsg{\Gamma}^{-1}(\ts{y}))
\det\nabla\Gamma^{-1}(\ts{y})\,\dif\ts{y}\\
&&+\int_{\Omega_\eps\setminus\Omega_\eps^\eta}\det\ts{A}\,\dif\ts{y},\\
&=&\int_{\Omega_\eps}\det\nabla\ts{u}\,\dif\ts{x}+\left(\det\ts{A}\right)
\left(\abs{\Omega}-\abs{\Omega^\eta}\right),\\
&=&\left(\det\ts{A}\right)\,\abs{\Omega}
-V_\eps+\left(\det\ts{A}\right)
\left(\abs{\Omega}-\abs{\Omega^\eta}\right),
\end{eqnarray*}
where
\[
\Omega^\eta=\set{\ts{x}\in\Omega\,:\,\mbox{dist}(\ts{x},\partial\Omega)>1-
\eta} .
\]
With $\eps$ fixed, choose $\eta(\eps)$ such that
\[
\left(\det\ts{A}\right)
\left(\abs{\Omega}-\abs{\Omega^{\eta(\eps)}}\right)=V_\eps-V.
\]
Note that $\eta(\eps)\nearrow1$ as $\eps\searrow0$. With this choice of $\eta(\eps)$, we get
that
\[
\int_{\Omega_\eps}\det\nabla\hat{\ts{u}}_{\eta(\eps)}\,\dif\ts{\ts{y}}=
\left(\det\ts{A}\right)\,\abs{\Omega}
-V.
\]
Combining this with result of part (i), we get that
$\hat{\ts{u}}_{\eta(\eps)}\in\mc{C}_{\ts{A}}^\eps$.
\end{enumerate}
Henceforth we set $\hat{\ts{u}}_\eps=\hat{\ts{u}}_{\eta(\eps)}$. Since
$\tsg{\Gamma}(\ts{x})=\ts{x}$ for $\ts{x}\in\partial\B_\eps(\ts{x}_0)$, we can
extend $\hat{\ts{u}}_\eps$ to $\Omega$ by setting
$\hat{\ts{u}}_\eps(\ts{x})=\ts{u}(\ts{x})$ for $\ts{x}\in B_\eps(\ts{x}_0)$ with
the resulting $\hat{\ts{u}}_\eps$ now in $W^{1,p}(\Omega)$ 
and with
$\left.\hat{\ts{u}}_\eps\right|_{\Omega_{\eps}}\in\mc{C}_{\ts{A}}^{\eps}$.

Now take $\set{\eps_j}$ with $\eps_j\To0$ and set $\hat{\ts{u}}_j=\hat{\ts{u}}_{\eps_j}$. Let
$\tsg{\Gamma}_j:\Omega_{\eps_{j}}\To\Omega_{\eps_{j}}^{\eta(\eps_j)}$ be the
corresponding diffeomorphism in the definition of $\hat{\ts{u}}_j$.
Since $\tsg{\Gamma}_j$ when restricted to $\partial\B_{\eps_{j}}(\ts{x}_0)$ is equal to the identity mapping, it follows that it can be extended continuously into $\B_{\eps_{j}}$ as the
identity. Thus $\nabla\tsg{\Gamma}_j^{-1}\To\ts{I}$ as $j\To\infty$ in $L^\infty(\Omega)$. These observations together with \eqref{H3} imply that
\[
W(\nabla\hat{\ts{u}}_{j})\le K[W(\nabla\ts{u})+1],
\]
a.e. in $\Omega$ for $j$ large enough. Now $\nabla\hat{\ts{u}}_j\To\nabla\ts{u}$ a.e. in $\Omega$. The above inequality and the Dominated Convergence
Theorem can be used now to conclude that
\[
\limite{j}{\infty}\int_{\Omega}W(\nabla\hat{\ts{u}}_{j}
(\ts{x}))\, \dif\ts{ x }
=\int_{\Omega}W(\nabla\ts{u}(\ts{x}))\,\dif\ts{x}.
\]
\end{prueba}
We now have the main result of this paper.
\begin{thm}\label{rcmconv}
Let the hypotheses in Lemma \ref{lem:3.1} hold. For
$V\in(0,(\omega_n/n)\det\ts{A})$, let $\set{\eps_j}$ be a sequence of
positive numbers converging to zero, and for each $\eps_j$, let $\ts{u}_j$ be
the corresponding minimizer given by Theorem \ref{thm:4.2} and satisfying
\eqref{constP}. Then $\set{\ts{u}_j}$ has a subsequence $\set{\ts{u}_{j_k}}$
such that for any $\delta>0$,
\[
\ts{u}_{j_k}\rightharpoonup\ts{u}_V\quad\mbox{in}\quad W^{1,p}(\Omega_\delta),
\]
where the function $\ts{u}_V$ is a solution of \eqref{AVmin0}, and with
\[
E(\ts{u}_V)=\lim_k\,E_{\eps_{j_k}}(\ts{u}_{j_k}).
\]
\end{thm}

\begin{prueba}
We let
\[
\mc{C}_{\ts{A}}^{\eps}\equiv\set{\ts{u}\in\mc{A}^\eps_{\ts{A}}\,|\,c_\eps(\ts{u}
)=0 }, \quad \eps\ge0.
\]
It follows from Lemma \ref{constsets} that these sets are non empty for $\eps$
small enough. Thus each $\ts{u}_j$ satisfies:
\[
E_{\eps_j}(\ts{u}_j)=
\min_{\ts{u}\in\mc{C}_{\ts{A}}^{\eps_j}}\int_{\Omega_{\eps_j}}
W(\nabla\ts{u}(\ts{x}))\,\dif\ts{x}=\min_{\ts{u}\in\mc{C}_{\ts{A}}^{\eps_j}}
E_{\eps_j}(\ts{u}).
\]
Now we fix an index $J\in\mathrm{N}$ and take $j>J$. It follows from hypothesis
(H2) on $W$ that for some constants $c_1>0$ and $c_2\in\Real$:
\[
E_{\eps_J}(\ts{u}_j)\ge
c_1\norm{\nabla\ts{u}_j}^p_{L^p(\Omega_{\eps_J})}+c_2,\quad
j>J.
\]
Again, it follows from (H2) that we may assume that $W$ is non negative. Hence
\[
E_{\eps_J}(\ts{u}_j)\le E_{\eps_j}(\ts{u}_j)\le C,\quad j>J,
\]
where the constant $C$ is given by Lemma \ref{constsets}. Combining this with
the inequality above, and since $\ts{u}_j=\ts{A}\ts{x}$ on
$\partial\Omega$,
we get that (for a subsequence) $\set{\ts{u}_{j}}$ converges
weakly in $W^{1,p}(\Omega_{\eps_J})$ to a function $\ts{u}^J$, and that
$\set{\det\nabla\ts{u}_{j}}$ converges weakly in
$L^1(\Omega_{\eps_J})$ to a function $\theta^J$. Since $p\in(n-1,n)$, it follows
from \cite[Theorem 4.2]{MuSp95}, that $\ts{u}^J$ satisfies condition INV,
$\theta^J=\det\nabla\ts{u}^J$, and $\det\nabla\ts{u}^J>0$ almost everywhere. By
choosing an appropriate diagonal sequence, it is shown in \cite{SiSpTi2006}
that there exists a subsequence $\set{\ts{u}_{j_k}}$ and a function $\ts{u}_V\in
W^{1,p}(\Omega)$ such that
\[
\ts{u}_{j_k}\rightharpoonup\ts{u}_V,\quad\mbox{in}\quad
W^{1,p}(\Omega_{\eps_J}).
\]
The results in \cite[Section 4.2]{SiSpTi2006} show that
$\ts{u}_V\in\mc{A}_\ts{A}$. 

It remains to show that $\ts{u}_V$ is a solution of \eqref{AVmin0}. By the
results quoted in the previous paragraph, we get that the
subsequence $\set{\ts{u}_{j_k}}$ has the property that
\[
\det\nabla\ts{u}_{j_k}\rightharpoonup\det\nabla\ts{u}_V,\quad\mbox{in}\quad
L^{1}(\Omega_{\eps_J}).
\]
Since $\ts{u}_{j_k}\in\mc{C}_{\ts{A}}^{\eps_{j_k}}$, we also have that
\[
\int_{\Omega_{\eps_{j_k}}}\det\nabla\ts{u}_{j_k}\,\dif\ts{x}=\left(\det\ts{A}
\right)\,\abs{\Omega}-V,\quad\forall k.
\]
Now we extend $\det\nabla\ts{u}_{j_k}$ to $\Omega$ as follows:
\[
g_k(\ts{x})=\left\{\begin{array}{ll}\det\nabla\ts{u}_{j_k}(\ts{x}),&\ts{x}\in
\Omega_{\eps_{j_k}},\\
0,&\ts{x}\in\Omega\setminus\Omega_{\eps_{j_k}}.\end{array}\right.
\]
Clearly $g_k\in L^1(\Omega)$ and since $\det\nabla\ts{u}_{j_k}>0$ a.e. in
$\Omega_{\eps_{j_k}}$, we get that
\[
\norm{g_k}_{L^1(\Omega)}=\int_{\Omega}
g_k\,\dif\ts{x}=\int_{\Omega_{\eps_{j_k}}}\det\nabla\ts{ u }_{j_k}\,\dif\ts{
x}=\left(\det\ts{A}\right)\,\abs{\Omega}-V,\quad\forall k .
\]
Writing
\[
\int_{\Omega}(\det\nabla\ts{u}_V-g_k)\,\dif\ts{x}=
\int_{\Omega_{\eps_J}}(\det\nabla\ts{u}_V-g_k)\,\dif\ts{x}+
\int_{\Omega\setminus\Omega_{\eps_J}}(\det\nabla\ts{u}_V-g_k)\,\dif\ts{x},
\]
we observe that the second term on the right is bounded by a constant times
$\abs{\Omega\setminus\Omega_{\eps_J}}$, and thus can be made arbitrarily small
by
taking $J$ sufficiently large. Once $J$ is fixed, the first term can be
make arbitrarily small for $k$ sufficiently large, as
$g_k=\det\nabla\ts{u}_{j_k}$ over $\Omega_J$ for $k$ sufficiently large
and by the weak convergence of $\set{\det\nabla\ts{u}_{j_k}}$ to
$\det\nabla\ts{u}_{V}$ in $L^1(\Omega_J)$. This essentially shows that
\[
\int_{\Omega}\det\nabla\ts{u}_V\,\dif\ts{x}=\limite{k}{\infty}\int_{\Omega}g_k\,
\dif\ts{x }=\left(\det\ts{A}\right)\,\abs{\Omega}-V,
\]
Hence $\ts{u}_V\in\mc{C}_{\ts{A}}^0$. We now show that $\ts{u}_V$ is a minimizer
over $\mc{C}_{\ts{A}}^0$.

For any $\ts{u}\in\mc{C}_{\ts{A}}^0$ and for the subsequence $\set{\eps_{j_k}}$
above, let $\set{\hat{\ts{u}}_{j_k}}$ be the corresponding sequence given by
Lemma \ref{lem:3.1} with the property that
\begin{equation}\label{limsq:1}
\limite{k}{\infty}\int_{\Omega}W(\nabla\hat{\ts{u}}_{j_k})
(\ts{x}))\, \dif\ts{ x }
=\int_{\Omega}W(\nabla\ts{u}(\ts{x}))\,\dif\ts{x}.
\end{equation}
As a function over $\Omega_{\eps_{j_k}}$, we have that
$\hat{\ts{u}}_{j_k}\in\mc{C}_{\ts{A}}^{\eps_{j_k}}$. Since $\ts{u}_{j_k}$ is
the minimizer over $\mc{C}_{\ts{A}}^{\eps_{j_k}}$, we have that
\begin{equation}\label{limsq:2}
\int_{\Omega_{\eps_{j_k}}}W(\nabla\ts{u}_{j_k}(\ts{x}))\,\dif\ts{x}\le
\int_{\Omega_{\eps_{j_k}}}W(\nabla\hat{\ts{u}}_{j_k}(\ts{x}))\,\dif\ts{x}.
\end{equation}
Let $N>0$ be given. For $k>N$ and the nonnegativity of $W$ we get that
\begin{equation}\label{limsq:4}
\int_{\Omega_{\eps_{j_N}}}W(\nabla\ts{u}_{j_k}(\ts{x}))\,\dif\ts{x}\le
\int_{\Omega_{\eps_{j_k}}}W(\nabla\ts{u}_{j_k}(\ts{x}))\,\dif\ts{x}.
\end{equation}
By the results in \cite{Ba82}, the functional $E_{\eps_{j_N}}(\cdot)$ (cf.
\eqref{regenergy}) is weakly lower semi--continuous over
$\mc{A}_{\ts{A}}^{\eps_{j_N}}$. Using this and since
$\ts{u}_{j_k}\rightharpoonup\ts{u}_V$ in $W^{1,p}(\Omega_{\eps_{j_N}})$,
we conclude that
\begin{equation}\label{limsq:3}
\int_{\Omega_{\eps_{j_N}}}W(\nabla\ts{u}_V(\ts{x}))\,\dif\ts{x}\le
\liminf_k \int_{\Omega_{\eps_{j_N}}}W(\nabla\ts{u}_{j_k}(\ts{x}))\,\dif\ts{x}.
\end{equation}
From the nonnegativity of $W$, it follows from \eqref{limsq:1}--\eqref{limsq:3}
that
\[
\int_{\Omega_{\eps_{j_N}}}
W(\nabla\ts{u}_V(\ts{x}))\,\dif\ts{x}\le\int_{\Omega}W(\nabla\ts{u}(\ts{x}))\,
\dif\ts {x}.
\]
Since $N$ is arbitrary, we can conclude that 
\[
\int_{\Omega}
W(\nabla\ts{u}_V(\ts{x}))\,\dif\ts{x}\le\int_{\Omega}W(\nabla\ts{u}(\ts{x}))\,
\dif\ts {x}.
\]
Since $\ts{u}\in\mc{C}_{\ts{A}}^0$ is arbitrary, we get that $\ts{u}_V$ is a
minimizer over $\mc{C}_{\ts{A}}^0$. If we set $\ts{u}=\ts{u}_V$ in
\eqref{limsq:1}, then we get as well that
\[
\int_\Omega
W(\nabla\ts{u}_V(\ts{x}))\,\dif\ts{x}=\liminf_k\int_{\Omega_{\eps_{j_k}}}
W(\nabla\ts { u } _ { j_k}(\ts{x}))\,\dif\ts{x},
\]
from which the result about the energies follows upon taking another
subsequence.
\end{prueba}

We now derive an expression for a weak form of the equilibrium equations for the
minimizer $\ts{u}_V$ in Theorem \ref{rcmconv}.

\begin{thm}\label{equivEQ}
Assume that \eqref{growthhyp2} and the hypotheses in Theorem \ref{rcmconv} hold.
Let $\ts{u}_V$ be the minimizer in Theorem \ref{rcmconv}. Then there exists
$\mu_V\in\Real$ such that
\begin{equation}\label{wfELV}
\int_{\Omega}\left[\td{W}{\ts{F}}(\nabla\ts{u}_V)
+\mu_V\left(\Adj\nabla\ts{u}_V\right)^T
\right]\cdot\nabla[\ts{v}(\ts{u}_V)]\,\dif\ts{x}=0,
\end{equation}
for all $\ts{v}\in C^1(\Real^n)$ with $\ts{v}=\ts{0}$ on
$\Real^n\setminus\mc{E}$, where $\mc{E}=\set{\ts{A}\ts{x}\,:\,\ts{x}\in
\Omega}$. Moreover, if $\ts{u}_V\in C^2(\Omega\setminus\set{\ts{x}_0})\cap
C^1(\overline{\Omega}\setminus\set{\ts{x}_0})$ with $\det\nabla\ts{u}_V>0$ in
$\Omega\setminus\set{\ts{x}_0}$, then
\begin{equation}\label{sfELV}
\Div\left[\displaystyle\td{W}{\ts{F}}(\nabla\ts{u}_V)
+\mu_V\left(\Adj\nabla\ts{u}_V\right)^T\right]=\ts{0},\quad\mbox{in~~}
\Omega\setminus\set{\ts{x}_0},
\end{equation}
and
\begin{equation}\label{sfELVibc}
\limite{\delta}{0}\,\int_{\partial\B_\delta(\ts{x}_0)}
\left(\left[\td{W}{\ts{F}}(\nabla\ts{u}_V)
+\mu_V(\Adj\nabla\ts{u}_V)^T
\right]\cdot\ts{n}\right)\cdot\ts{v}(\ts{u}_V)\,\dif s(\ts{x})=0.
\end{equation}
\end{thm}

\begin{prueba}
Let $\set{\ts{u}_{j_k}}$ be the subsequence given by Theorem \ref{rcmconv} such
that for any $\delta>0$,
\begin{subequations}\label{convgsubseq}
\begin{eqnarray}
\ts{u}_{j_k}&\rightharpoonup&\ts{u}_V\quad\mbox{in}\quad
W^{1,p}(\Omega_\delta),\label{wconvgu}\\
\det\nabla\ts{u}_{j_k}&\rightharpoonup&\det\nabla\ts{u}_V\quad\mbox{in}\quad
L^{1}(\Omega_\delta),\label{wcongdet}
\end{eqnarray}
\end{subequations}
and with
\begin{equation}\label{convgener}
E(\ts{u}_V)=\lim_k\,E_{\eps_{j_k}}(\ts{u}_{j_k}).
\end{equation}
(We recall that $\ts{u}_{j_k}$ is the corresponding minimizer given by
Theorem \ref{thm:4.2} corresponding to $\eps_{j_k}$.) From Theorem
\ref{equivEQeps} we have that
\[
\int_{\Omega_{\eps_{j_k}}}\left[\td{W}{\ts{F}}(\nabla\ts{u}_{\eps_{j_k}})
\nabla\ts{u}_{\eps_{j_k}}^T+\mu_{\eps_{j_k}}(\det\nabla\ts{u}_{\eps_{j_k}})
\ts{I}\right]\cdot\nabla\ts{v}(\ts{u}_{\eps_{j_k}})\,\dif\ts{x}=0.
\]
Let $\ts{g}_{\eps_{j_k}}$ and $\ts{h}_{\eps_{j_k}}$ represent the extensions (by
the zero tensor) over $\Omega$ of the functions
\[
\td{W}{\ts{F}}(\nabla\ts{u}_{\eps_{j_k}})
\nabla\ts{u}_{\eps_{j_k}}^T,\quad (\det\nabla\ts{u}_{\eps_{j_k}})
\ts{I},
\]
respectively. The sequence $\set{\ts{u}_{\eps_{j_k}}}$ can be extended (cf.
\cite[p. 746]{SiSpTi2006}) as well to a sequence
$\set{\tilde{\ts{u}}_{\eps_{j_k}}}$ over $\Omega$ in such
way that \eqref{convgsubseq} holds now over $\Omega$, with
$\tilde{\ts{u}}_{\eps_{j_k}}= \ts{u}_{\eps_{j_k}}$ over $\Omega_{\eps_{j_k}}$,
and
with $\tilde{\ts{u}}_{\eps_{j_k}}(\Omega)\subset\mc{E}$. Using
\eqref{growthhyp2},
\eqref{convgsubseq}. and \eqref{convgener} together with the Dominated
convergence Theorem, we get that
\begin{eqnarray*}
\int_{\Omega}\td{W}{\ts{F}}(\nabla\ts{u}_V)\nabla\ts{u}_V^T
\cdot\nabla\ts{v}(\ts{u}_V)\,\dif\ts{x}&=&
\limite{k}{\infty}\int_{\Omega}\ts{g}_{\eps_{j_k}}
\cdot\nabla\ts{v}(\tilde{\ts{u}}_{\eps_{j_k}})\,\dif\ts{x},\\
\int_{\Omega}(\det\nabla\ts{u}_V)\ts{I}
\cdot\nabla\ts{v}(\ts{u}_V)\,\dif\ts{x}&=&
\limite{k}{\infty}\int_{\Omega}\ts{h}_{\eps_{j_k}}
\cdot\nabla\ts{v}(\tilde{\ts{u}}_{\eps_{j_k}})\,\dif\ts{x}.
\end{eqnarray*}
But
\begin{eqnarray*}
&0=\displaystyle\int_{\Omega_{\eps_{j_k}}}\left[\td{W}{\ts{F}}
(\nabla\ts{u}_{\eps_{j_k } } )
\nabla\ts{u}_{\eps_{j_k}}^T+\mu_{\eps_{j_k}}(\det\nabla\ts{u}_{\eps_{j_k}})
\ts{I}\right]\cdot\nabla\ts{v}(\ts{u}_{\eps_{j_k}})\,\dif\ts{x}&\\
&~~~~~~~~~~~~~~~~=\displaystyle\int_{\Omega}\left[\ts{g}_{\eps_{j_k}}+\mu_{\eps_
{j_k
}}\ts{h}_{ \eps_ { j_k } } \right ]
\cdot\nabla\ts{v}(\tilde{\ts{u}}_{\eps_{j_k}})\,\dif\ts{x},&
\end{eqnarray*}
for all $k$. Combining this with the two limits above we that there exists
$\mu_V\in\Real$ such that $\mu_{\eps_{j_k}}\To\mu_V$ and with
\[
\int_{\Omega}\left[\td{W}{\ts{F}}(\nabla\ts{u}_V)\nabla\ts{u}_V^T+
\mu_V(\det\nabla\ts{u}_V)\ts{I}\right]
\cdot\nabla\ts{v}(\ts{u}_V)\,\dif\ts{x}=0,
\]
from which \eqref{wfELV} follows. 

Now assume that $\ts{u}_V\in C^2(\Omega\setminus\set{\ts{x}_0})\cap
C^1(\overline{\Omega}\setminus\set{\ts{x}_0})$ with $\det\nabla\ts{u}_V>0$ in
$\Omega\setminus\set{\ts{x}_0}$. The proof of \eqref{sfELV} is similar to the
one given in \cite[Theorem 5.1]{SiSp2000a} and thus we omit it.
Let $\delta>0$ be given. If we multiply \eqref{sfELV} by $\ts{v}(\ts{u}_V)$,
where $\ts{v}\in C^1(\Real^n)$ with $\ts{v}=\ts{0}$ on
$\Real^n\setminus\mc{E}$, and integrate by parts over $\Omega_\delta$, we get
that
\begin{eqnarray*}
&\displaystyle\int_{\Omega_\delta}\left[\td{W}{\ts{F}}(\nabla\ts{u}_V)
+\mu_V\left(\Adj\nabla\ts{u}_V\right)^T
\right]\cdot\nabla[\ts{v}(\ts{u}_V)]\,\dif\ts{x}\quad\quad\quad\quad&\\
&\quad\quad=\displaystyle\int_{\partial\B_\delta(\ts{x}_0)}
\left(\left[\td{W}{\ts{F}}(\nabla\ts{u}_V)
+\mu_V(\Adj\nabla\ts{u}_V)^T
\right]\cdot\ts{n}\right)\cdot\ts{v}(\ts{u}_V)\,\dif s(\ts{x}).&
\end{eqnarray*}
Taking the limit as $\delta\searrow0$ and using \eqref{wfELV} we get that
\eqref{sfELVibc} holds.

\end{prueba}

\begin{remark}
One could try to prove Theorem \ref{equivEQ} by using variations directly onto
the functional in \eqref{AVmin0}. However this approach would require
constructing variations that preserve the volume constraint $c_0(\ts{u})=0$
(cf. \eqref{const0}). Our proof using Theorem \ref{equivEQeps} avoids the
technical complications of constructing such variations.
\end{remark}

\begin{remark}
In terms of the Cauchy stress tensor:
\[
\ts{T}(\ts{u}_V)=(\det\nabla\ts{u}_V)^{-1}\,\td{W}{\ts{F}}(\nabla\ts{u}
_V)(\nabla\ts{u}_V)^T,
\]
\eqref{sfELVibc} is equivalent to:
\[
\limite{\delta}{0}\,
\int_{\ts{u}_V(\partial\B_\delta(\ts{x}_0))}\left(\left[\ts{T}(\ts{y})+\mu_V\ts{
I }
\right]\cdot\tilde{\ts{n}}(\ts{y})\right)\cdot\ts{v}(\ts{y})\,\dif s(\ts{y})=0,
\]
where $\tilde{\ts{n}}$ is the unit normal to
$\ts{u}_V(\partial\B_\delta(\ts{x}_0))$. If $H$ is the region of volume $V$
occupied by the cavity induced by $\ts{u}_V$, then the limit above can be
replaced by the corresponding integral over $\partial H$. It follows now that
\[
\int_{\partial H}\left(\left[\ts{T}(\ts{y})+\mu_V\ts{
I }
\right]\cdot\tilde{\ts{n}}(\ts{y})\right)\cdot\ts{v}(\ts{y})\,\dif s(\ts{y})=0,
\]
for all $\ts{v}\in C^1(\Real^n)$. Thus
\begin{equation}\label{tracemult}
\left[\ts{T}(\ts{y})+\mu_V\ts{I}\right]\cdot\tilde{\ts{n}}(\ts{y})=0, \mbox{ 
over }\partial H,
\end{equation}
in the sense of trace.
\end{remark}
\section{Numerical results}\label{numres}
In this section we describe some of the elements of a numerical procedure, 
based on the results of the previous sections, to compute a minimizer of
\eqref{AVmin0}. In addition we work a numerical example in which we check the
convergence as $\eps\searrow0$ predicted by Theorem \ref{rcmconv}.

For given values of $\eps, V$, we use the method outlined in Theorem
\ref{thm:4.2}
to compute the minimizer $\ts{u}_{\eps}$ in \eqref{constP}. The minimizers in
\eqref{penalE} (dropping the subscript ``$j$'') are computed using the
\textit{gradient flow
equation}\footnote{It follows from \eqref{aux121a} that the Euler--Lagrange
equations for \eqref{penalE} are
formally given by equating to zero the right hand side of \eqref{GFequ}.}:
\begin{equation}\label{GFequ}
\Delta\ts{u}_t=-\Div\left[\displaystyle\td{W}{\ts{F}}(\nabla\ts{u})+
(\mu+\eta c_\eps(\ts{u}))(\Adj\nabla\ts{u})^t\right],\mbox{ in~~}\Omega_\eps,
\end{equation}
where for all $t\ge0$, $\ts{u}(\ts{x},t)=\ts{A}\ts{x}$ over $\partial\Omega$ and
\begin{equation}\label{GFequbc}
\displaystyle\left[\nabla\ts{u}_t+\td{W}{\ts{F}}(\nabla\ts{u})+
(\mu+\eta c_\eps(\ts{u}))(\Adj\nabla\ts{u})^t\right]\cdot \ts{n}=\ts{0},\mbox{
~on~}\partial\B_\eps(\ts{x}_0).
\end{equation}
The gradient flow equation leads to a descent method for
the solution of \eqref{penalE}. (For more
details about the gradient flow method and its properties we refer to
\cite{Neu1997}, and for its use in problems
leading to cavitation see \cite{HX}.) After discretization of the
partial derivative with respect to ``$t$", one can use a
finite element method to solve the resulting flow
equation. In particular, if we let $\Delta t>0$
be given, and set $t_{i+1}=t_i+\Delta t$ where $t_0=0$, we can approximate
$\ts{u}_t(\ts{x},t_i)$ with:
\[
\ts{z}_i(\ts{x})= \dfrac{\ts{u}_{i+1}(\ts{x})-\ts{u}_i(\ts{x})}{\Delta
t},
\]
where $\ts{u}_i(\ts{x})=\ts{u}(\ts{x},t_i)$, etc.. (We take $\ts{u}_0(\ts{x})$
to be some initial deformation satisfying the boundary condition on $\partial
\Omega$, e.g., $\ts{A}\ts{x}$.) Inserting this approximation into the weak form
of 
\eqref{GFequ}, \eqref{GFequbc}, and evaluating the right hand side of
\eqref{GFequ} at $\ts{u}=\ts{u}_i$, we arrive at the following iterative
formula:
\begin{equation}\label{floweq1}
\int_{\Omega_\eps}\nabla\ts{z}_i:\nabla\ts{v}\,\dif\ts{x}+
\int_{\Omega_\eps}\bigg[\td{W}{\ts{F}}(\nabla\ts{u}_i)
+(\mu+\eta
c_\eps(\ts{u_i}))(\Adj\nabla\ts{u}_i)^t\bigg]:\nabla\ts{v}\,\dif\ts{x}=0
,
\end{equation}
for all $\ts{v}$ vanishing on $\partial\Omega$ and sufficiently smooth so
that the
integrals above are well defined. Given $\ts{u}_i$, one can solve the above
equation for $\ts{z}_i$ via some
finite element scheme, and then set $\ts{u}_{i+1}=\ts{u}_i+\Delta
t\,\ts{z}_i$. This process is repeated for $i=0,1,\ldots$, until
$\ts{u}_{i+1}-\ts{u}_i$ is ``small'' enough ($10^{-3}$ in the calculations
below), or some maximum value of
``$t$'' is reached, declaring the last $\ts{u}_i$ as an approximation of
$\ts{u}_{\eps}$. This whole process is repeated for smaller values of
of $\eps$, to obtain as a result an approximation of the minimizer $\ts{u}_V$ in
\eqref{AVmin0}.

For the computations we used the stored energy function \eqref{sefexample} in
which:
\[
h(d)=c_1d^{e_1}+c_2d^{ -e_2},
\]
where $c_1,c_2\ge0$ and $e_1,e_2>0$. The reference configuration
is stress free provided:
\[
c_2=\dfrac{\mu(\sqrt{n})^{q-2}+c_1e_1}{e_2}.
\]
The case $\mu=0$ in \eqref{sefexample} is called an \textit{elastic fluid}. 

For an elastic fluid in which $\Omega=\B\equiv \B_1(\ts{0})$ and
$\ts{x}_0=\ts{0}$, the minimizer
$\ts{u}_V$ in \eqref{AVmin0} is given by (see \cite{NeSi2011b}):
\[
\ts{u}_V(\ts{x})=[dR^n +(1-d)]^{1/n}\, \dfrac{\ts{A}\ts{x}}{R},\quad
R=\norm{\ts{x}},
\]
where $d$ is given by
\[
d=1-\dfrac{nV}{\omega_n\det\ts{A}}.
\]
($V$ is assumed to be sufficiently small as to guarantee that $d>0$.) It
follows that $\det\nabla\ts{u}_V=d\det\ts{A}$. Thus we have that
\[
E(\ts{u}_V)=\int_{B}h(\det\nabla\ts{u}_V)\,\dif\ts{x}=\frac{\omega_n}{n}\,
h(d\det\ts{A}),
\]
where $\omega_n$ denotes the area of the unit sphere in $\Real^n$.
We now consider the particular case in which $n=2$, $c_1=1$, $e_1=2$,
$e_2=1$, $V=\pi(0.15)^2$, and $\ts{A}=\mbox{diag}(1.1,1.4)$. Using the formulas
above, we get that
\[
E(\ts{u}_V)=\pi\,h((1.1)(1.4)-0.15^2)=11.3750.
\]
For the parameters in Theorem \ref{thm:4.2} we used $\gamma=0.25$, $\beta=2$,
with the stopping criteria in \eqref{penalITER} given by
$|\mu_{j+1}-\mu_{j}|<10^{-3}|\mu_j|$, and for the solution of the sub--problems
\eqref{floweq1} we used the package freefem++ (see \cite{He2012}). We show in
Table \ref{tab:1} the results
in this case for the method described at the beginning of this section. For
each $\eps$ we clearly see the penalty--multiplier iterations converging as
predicted in Theorem \ref{thm:4.2}. Note that the penalty parameters do not
become too large, thus avoiding the ill--conditioning associated with large
values of these parameters. As we move down along the table and look at the last
computed energy for each $\eps$, we see that these values are approaching the
exact energy $11.3750$, to within the convergence tolerances in the gradient
flow and penalty multiplier iterations and finite element approximation, in
accordance with the result in Theorem \ref{rcmconv}.
\begin{table}
\scriptsize{\begin{center}
\begin{tabular}{|c|c|c|c|c|c|}\hline
$\eps$&$j$&$c_{\eps}(\ts{u}_j)$&$E_{\eps,\mu_j,\eta_j}(\ts{u}_j)$&  
$\mu_j$&$\eta_j$\\\hline\hline
0.1&0&   -0.435566&   10.5824&   0&     5\\\cline{2-6}
&1&   -0.108628&   11.3009&   -2.17783&    10 \\\cline{2-6}
&2&   -0.0182433&   11.3627&   -3.26411&     10  \\\cline{2-6}
&3&   -0.00175846&   11.3637&   -3.44654&   10  \\\cline{2-6}
&4&   0.00123636&   11.3636&   -3.46413&  10 \\\cline{2-6}
&5&   0.00119778&   11.3636&   -3.45177&  20\\\cline{2-6}
&6&   5.80566e-05&   11.3636&   -3.42781&  40\\\hline\hline
0.05&0&   -0.408565&   10.6179&   0&   5\\\cline{2-6}
&1&   -0.114678&   11.2988&   -2.04282&  10\\\cline{2-6}
&2&   -0.0116487&   11.3693&   -3.18961&  20\\\cline{2-6}
&3&   0.000242224&   11.3699&   -3.42258&  20\\\cline{2-6}
&4&   0.00100667&   11.3699&   -3.41774&  20\\\cline{2-6}
&5&   0.000804023&   11.3699&   -3.3976&  40\\\cline{2-6}
&6&   -3.24128e-05&   11.3699&   -3.36544&  80\\\hline\hline
0.025&0&   -0.201328&   10.8508&   0&   5\\\cline{2-6}
&1&   -0.193863&   11.149&   -1.00664&  10\\\cline{2-6}
&2&   -0.0446432&   11.3758&   -2.94527&  20\\\cline{2-6}
&3&   0.00475689&   11.3697&   -3.83813&  20\\\cline{2-6}
&4&   0.015619&   11.3684&   -3.74299&  20\\\cline{2-6}
&5&   0.00201459&   11.3716&   -3.43061&  40\\\cline{2-6}
&6&   -9.22278e-05&   11.3717&   -3.35003& 40\\\cline{2-6}
&7&   -9.20066e-05&   11.3717&   -3.35372&  40\\\cline{2-6}
&8&   -7.96496e-05&   11.3717&   -3.3574&  80\\\cline{2-6}
&9&   -5.49799e-06&   11.3717&   -3.36377&  160\\\hline\hline
0.0125&0&   -0.0427461&   11.2392&   0&   5\\\cline{2-6}
&1&   -0.0916523&   11.1388&   -0.213731&  10\\\cline{2-6}
&2&   -0.0956048&   11.2616&   -1.13025&  20\\\cline{2-6}
&3&   -0.0281877&   11.3806&   -3.04235&  40\\\cline{2-6}
&4&   0.00372948&   11.3698&   -4.16986&  80\\\cline{2-6}
&5&   0.00570313&   11.3705&   -3.8715&  80\\\cline{2-6}
&6&   0.000615149&   11.3721&   -3.41525&  160\\\cline{2-6}
&7&   -0.000212807&   11.3721&   -3.31683&  160\\\cline{2-6}
&8&   -3.51888e-05&   11.3721&   -3.35088&  320\\\cline{2-6}
&9&   1.10978e-06&   11.3721&   -3.36214&  320\\\hline\hline
0.00625&0&   0.0374087&   11.503&   0&   5\\\cline{2-6}
&1&   0.0174181&   11.4372&   0.187044&  10\\\cline{2-6}
&2&   0.00767417&   11.4034&   0.361225&  20\\\cline{2-6}
&3&   0.00229717&   11.3836&   0.514708&  40\\\cline{2-6}
&4&   -0.00539297&   11.3545&   0.606595&  80\\\cline{2-6}
&5&   -0.0106619&   11.3448&   0.175157&  160\\\cline{2-6}
&6&   -0.00686916&   11.3674&   -1.53074&  320\\\cline{2-6}
&7&   0.000534898&   11.3721&   -3.72887&  640\\\cline{2-6}
&8&   5.6007e-05&   11.3723&   -3.38654&  640\\\cline{2-6}
&9&   5.51664e-07&   11.3723&   -3.35069&  640\\\hline
\end{tabular} 
\end{center}
\caption{Convergence of the regularized penalty--multiplier minimizers for the
case of a two dimensional elastic fluid.}\label{tab:1}}
\end{table}

\section{Final Comments}
In \cite[Proposition 6.1]{NeSi2011b} the authors introduced a regularized
penalty method for approximating solutions of \eqref{AVmin0}. They anticipated
without proof, the convergence of the corresponding regularised minimizers to a
solution of \eqref{AVmin0}. The result in Theorem \ref{rcmconv} fills that gap
and by using a penalty-multiplier scheme we also obtained information on the
Lagrange multiplier corresponding to the integral constraint on the volume
of the hole produced in the deformed configuration. Moreover, the use of the
penalty--multiplier technique leads to a more a stable numerical scheme as
compared to a standard penalty method, as in general one achieves convergence to
a minimum without having to make the penalty parameter excessively large, which
could lead to numerical ill conditioning. We anticipate that the scheme
introduced in this paper can be used as part of an efficient numerical scheme
for the computation of the volume derivative introduced in \cite{NeSi2011b}.

\appendix
\section{The integral constraint}\label{regProb}
In this section of the appendix we give a characterization of the condition on
the distributional determinant in \eqref{admisset} in terms of an integral
constraint.
\begin{prop}\label{Detconst}
For any $\bu\in\mc{A}_{\mathbf{A},V}$, the following integral constraint holds:
\begin{equation}\label{intconst}
\int_{\Omega}\det\nabla\ts{u}\,\dif\ts{x}=\left(\det\ts{A}\right)\,\abs{\Omega}
-V.
\end{equation}
Moreover, for any $\ts{u}\in C^2(\Omega\setminus\set{\ts{x}_0})\cap
C^1(\overline{\Omega}\setminus\set{\ts{x}_0})$ that satisfies \eqref{intconst},
with $\det\nabla \mathbf{u}>0$ a.e., $\bu(\bx )=\mathbf{A}\bx$ on
$\partial \Omega$, and $\bu$ satisfies INV on $\Omega$, it follows that 
\begin{equation}\label{eqA2}
\mathrm{Det}\nabla \bu = \mathrm{det}\nabla \bu +V \delta _{\ts{x}_0}.
\end{equation}
\end{prop}
\begin{prueba}
The first part of the proposition follows essentially from \cite[Lemma
8.1]{MuSp95} but we refer to \cite{Va2007} for a proof as well.

For the second part of the proposition, let $\ts{u}\in
C^2(\Omega\setminus\set{\ts{x}_0})\cap
C^1(\overline{\Omega}\setminus\set{\ts{x}_0})$ with $\det\nabla \mathbf{u}>0$
a.e., $\bu(\bx )=\mathbf{A}\bx$ on
$\partial \Omega$, and $\bu$ satisfies INV on $\Omega$. Assume that
$\ts{u}$ satisfy \eqref{intconst}. For any $\phi\in C_0^\infty(\Omega)$ and
$\delta>0$ sufficiently small,
\begin{eqnarray*}
\int_{\Omega\setminus
\overline{\B_{\delta}(\ts{x}_0)}}\nabla\phi\cdot(\Adj\nabla\ts{u})\ts{u}\,\dif\ts
{ x } &=&\int_
{\partial({\Omega\setminus
\overline{\B_{\delta}(\ts{x}_0)}})}\phi\,((\Adj\nabla\ts{u})\ts{u})\cdot\ts{n}\,
\dif s\\&&-
\int_{\Omega\setminus
\overline{\B_{\delta}(\ts{x}_0)}}\phi\,\Div((\Adj\nabla\ts{u})\ts{u})\,\dif\ts{x}
, \\
&=&-\phi(\ts{x}_\delta)\int_{\partial
\B_{\delta}(\ts{x}_0)}((\Adj\nabla\ts{u})\ts{u})\cdot\ts{n}\,\dif s\\&&-
n\int_{\Omega\setminus\overline{
\B_{\delta}(\ts{x}_0)}}\phi\,\det\nabla\ts{u}\,\dif\ts{x},
\end{eqnarray*}
where $\ts{x}_\delta\in\partial B_{\delta}(\ts{x}_0)$ and $\ts{n}$ in the last
boundary integral is the unit outer normal to $\partial B_{\delta}(\ts{x}_0)$.
Thus:
\begin{eqnarray}
<\mbox{Det}\nabla
\bu,\phi>&=&\limite{\delta}{0}\left[-\dfrac{1}{n}\int_{\Omega\setminus
\overline{\B_{\delta}(\ts{x}_0)}}\nabla\phi\cdot(\Adj\nabla\ts{u})\ts{u}\,\dif\ts
{ x}
\right],\nonumber\\
&=&\phi(\ts{x}_0)\limite{\delta}{0}\left[\dfrac{1}{n}\int_{\partial
\B_{\delta}(\ts{x}_0)}((\Adj\nabla\ts{u})\ts{u})\cdot\ts{n}\,\dif
s\right]\nonumber\\&&+
\int_{\Omega}\phi\,\det\nabla\ts{u}\,\dif\ts{x}.\label{eqA1}
\end{eqnarray}
Note that:
\begin{eqnarray}
\int_{\Omega\setminus\overline{
B_{\delta}(\ts{x}_0)}}\det\nabla\ts{u}\,\dif\ts{x}&=&
\dfrac{1}{n}\int_{\partial({\Omega\setminus\overline{ B_{\delta}(\ts{x}_0)}})}
((\Adj\nabla\ts{u})\ts{u})\cdot\ts{n}\,\dif s,\nonumber\\
&=&\dfrac{1}{n}\int_{\partial \Omega}((\Adj\nabla\ts{u})\ts{u})\cdot\ts{n}\,\dif
s\nonumber\\
&&-\dfrac{1}{n}\int_{\partial B_{\delta}(\ts{x}_0)}
((\Adj\nabla\ts{u})\ts{u})\cdot\ts{n}\,\dif s.\label{eqA3}
\end{eqnarray}
Assuming that \eqref{intconst} holds, we get that:
\begin{eqnarray*}
\limite{\delta}{0}\left[\dfrac{1}{n}\int_{\partial
B_{\delta}(\ts{x}_0)}((\Adj\nabla\ts{u})\ts{u})\cdot\ts{n}\,\dif
s\right]&=&\dfrac{1}{n}\int_{\partial\Omega}((\Adj\nabla\ts{u})\ts{u})
\cdot\ts{n}\,\dif s
-\int_{\Omega}\det\nabla\ts{u}\,\dif\ts{x},\\
&=&\dfrac{1}{n}\int_{\partial\Omega}
((\Adj\nabla\ts{u})\ts{A}\ts{x})\cdot\ts{n}\,\dif
s+V-\left(\det\ts{A}\right)\,\abs{\Omega}.
\end{eqnarray*}
Since $\ts{u}\in C^1(\overline{\Omega}\setminus\set{\ts{x}_0})$ and
$\ts{u}(\ts{x})=\ts{A}\ts{x}$ for all
$\ts{x}\in\partial\Omega$, we have that $\nabla\ts{u}(\ts{x})=\ts{A}$ for all
$\ts{x}\in\partial\Omega$. Hence
\begin{eqnarray}
\int_{\partial
\Omega}((\Adj\nabla\ts{u})\ts{A}\ts{x})\cdot\ts{n}\,\dif s&=&
\int_{\partial\Omega} ((\Adj\ts{A})\ts{A}\ts{x})\cdot\ts{n}\,\dif s,\nonumber\\
&=&\det\ts{A}\int_{\partial\Omega}\ts{x}\cdot\ts{n}\,\dif
s=n\left(\det\ts{A}\right)\abs{\Omega}.\label{eqA4}
\end{eqnarray}
Using this above we get that
\[
\limite{\delta}{0}\left[\dfrac{1}{n}\int_{\partial
B_{\delta}(\ts{x}_0)}((\Adj\nabla\ts{u})\ts{u})\cdot\ts{n}\,\dif
s\right]=V.
\]
Thus combining this with \eqref{eqA1} we get that:
\[
<\mbox{Det}\nabla
\bu,\phi>=\int_{\Omega}\phi\,\det\nabla\ts{u}\,\dif\ts{x}+V\phi(\ts{x}_0),
\quad\forall\,\,\phi\in C_0^\infty(\Omega),
\]
i.e., that \eqref{eqA2} holds.
\end{prueba}

\section{The constrained admissible set}\label{constAS}
We now show that for $\eps$ sufficiently small, the admissible
sets in \eqref{AVmineps} are non--empty. In the proof we make use of the following result. let $\B_1(\ts{0})$ be the unit ball with centre at the origin and for any $d\in(0,1)$, define
\begin{equation}\label{basic_cav}
 \ts{u}_d(\ts{x})=[d R^n+(1-d)]^{1/n}\,\dfrac{\ts{x}}{R},
\quad R=\norm{\ts{x}},\quad\ts{x}\in\B_1(\ts{0}).
\end{equation}
It follows from \cite[Lemma 4.1]{Ba82}, that $\ts{u}_d\in W^{1,p}(\B_1(\ts{0}))$ for $p\in[1,n)$. An easy calculation shows as well that $\det\nabla\ts{u}_d=d$.
\begin{lem}\label{constsets}
There exists $V_0, \eps_0>0$ such that
\[
\mc{C}_{\ts{A}}^{\eps}\equiv\set{\ts{u}\in\mc{A}^\eps_{\ts{A}}\,|\,c_\eps(\ts{u}
)=0 } \ne\emptyset ,
\]
for all $\eps\in[0,\eps_0)$ and $0<V<V_0$. Moreover, if $W$ is nonnegative
and for any $0<\gamma<\delta$ there exists a
constant $K>0$ such that
\begin{equation}\label{hipW1}
W(\ts{F})\le K(\norm{\ts{F}}^p+1),\quad\mbox{whenever  }
\det\ts{F}\in[\gamma,\delta],
\end{equation}
then for any nonnegative sequence $\eps_j\To0$,
there exists a sequence $\ts{z}_{\eps_{j}}\in
\mc{C}_{\ts{A}}^{\eps_{j}}$ such that
\[
E_{\eps_{j}}(\ts{z}_{\eps_{j}})\le C,\quad \forall\,\,j,
\]
for some constant $C>0$.
\end{lem}
\begin{prueba}
Let $\eta>0$ be such that $\B_\eta(\ts{x}_0)\subset\Omega$ and define 
\[
V_0=\eta^n\,\dfrac{\omega_n}{n}\,\det\ts{A},\quad
\eps^n_0=\dfrac{nV}{\omega_n\det\ts{A}}.
\]
For $0<V<V_0$ and $\eps\in[0,\eps_0)$, let
\[
d_\eps=\dfrac{1}{\eta^n-\eps^n}\left(\eta^n-\dfrac{nV}{\omega_n\det\ts{A}}
\right).
\]
It follows now that $d_\eps\in(0,1)$ and that $B_\eps(\ts{x}_0)\subset
\B_\eta(\ts{x}_0)$. Now define $\ts{w}_\eps=\ts{u}_{d_\eps}$ and let
\begin{eqnarray*}
\ts{v}_\eps(\ts{x})=\left\{\begin{array}{rcl}
\eta \ts{A}\ts{w}_\eps\left(\frac{\ts{x}-\ts{x}_0}{\eta}+\ts{x}_0\right)
+\ts{A}\ts{x}_0-\eta\ts{A}\ts{x}_0&,&\ts{x}\in\B_\eta(\ts{x}_0),\\
\ts{A}\ts{x}&,&\ts{x}\in\Omega\setminus\B_\eta(\ts{x}_0).
\end{array}\right.
\end{eqnarray*}
Clearly $\ts{v}_\eps(\ts{x})=\ts{A}\ts{x}$ on $\partial\Omega\cup\partial\B_\eta(\ts{x}_0)$. Since $\ts{w}_\eps\in W^{1,p}(\B_1(\ts{0}))$ for $p\in[1,n)$, it follows that $\ts{v}_\eps\in W^{1,p}(\Omega)$ for $p\in[1,n)$. Hence
in particular, if $\ts{z}_\eps=\left.\ts{v}_\eps\right|_{\Omega_\eps}$, then
$\ts{z}_\eps\in\mc{A}^\eps_{\ts{A}}$ for $\eps\ge0$. That
$\ts{z}_0\in\mc{A}^0_{\ts{A}}$ follows from the following result:
\[
\mathrm{Det}\nabla\ts{w}_0=\mathrm{det}\nabla\ts{w}_0
+\dfrac{\omega_n}{n}(1-d_0)\delta_{\ts{0}},
\]
where $\delta_{\ts{0}}$ is the Dirac delta measure at the origin. An easy calculation now shows that
\begin{eqnarray*}
\det\nabla\ts{z}_{\eps}=\left\{\begin{array}{rcl}
d_\eps\det\ts{A}&,&\ts{x}\in\B_\eta(\ts{x}_0)\setminus\B_\eps(\ts{x}_0),\\
\det\ts{A}&,&\ts{x}\in\Omega\setminus\B_\eta(\ts{x}_0).
\end{array}\right.
\end{eqnarray*}
Using the definition of $d_\eps$, we get now that 
\begin{eqnarray*}
\int_{\Omega_\eps}\det\nabla\ts{z}_{\eps}\,\dif\ts{x}&=&
\int_{\B_\eta(\ts{x}_0)\setminus\B_\eps(\ts{x}_0)}\det\nabla\ts{z}_{\eps}\,
\dif\ts{x}+
\int_{\Omega\setminus\B_\eta(\ts{x}_0)}\det\nabla\ts{z}_{\eps}\,\dif\ts{x},\\
&=&\left(\det\ts{A}\right)\dfrac{\omega_n}{n} d_\eps(\eta^n-\eps^n)+
\det\ts{A}\left(\abs{\Omega}-\dfrac{\omega_n}{n}\eta^n\right),\\
&=&\left(\det\ts{A}\right)\abs{\Omega}-V,
\end{eqnarray*}
that is, $c_\eps(\ts{z}_{\eps})=0$. Hence
$\ts{z}_\eps\in\mc{C}_{\ts{A}}^{\eps}$.

For the second part of the lemma, we observe that for any nonnegative sequence
$(\eps_j)$ converging to zero, we can conclude from \cite[Eqn. (4.5)]{Ba82} that
$(\ts{v}_{\eps_j})$ is a Cauchy sequence in $W^{1,p}(\Omega)$. Since
$\ts{v}_{\eps_j}\To\ts{v}_0$ a.e., we have that $\ts{v}_{\eps_j}\To\ts{v}_0$ in
$W^{1,p}(\Omega)$. From \eqref{hipW1} and since
$(\det\nabla\ts{v}_{\eps_j})$ converges a.e. to
\[
\left\{\begin{array}{rcl}
d_0\det\ts{A}&,&\ts{x}\in\B_\eta(\ts{x}_0),\\
\det\ts{A}&,&\ts{x}\in\Omega\setminus\B_\eta(\ts{x}_0).
\end{array}\right.
\]
we get that
\[
E_{\eps_j}(\ts{z}_{\eps_j})=\int_{\Omega_{\eps_j}}W(\nabla\ts{z}_{\eps_j})\,
\dif\ts {
x }
\le \int_{\Omega}W(\nabla\ts{v}_{\eps_j})\,\dif\ts{x}\le
K(\norm{\nabla\ts{v}_{\eps_j}}^p+1)\le C,
\]
for some constant $C>0$.
\end{prueba}

\noindent\textbf{Acknowledgements:}
This research was sponsored in part by an International Joint Project
Grant from the Royal Society of London. The work of Negr\'on--Marrero also was
sponsored in part by the NSF--PREM Program of the UPRH (Grant No.
DMR--1523463).

\end{document}